\documentclass[12pt]{article} \baselineskip=.2cm \textwidth=165mm \textheight=22cm \voffset
-1.5cm \hoffset -1cm

\usepackage{amssymb}
\usepackage{amsfonts}
\usepackage{latexsym}
\usepackage{amsthm}

\usepackage[english]{babel}

\newcommand{\er}[1]{{\rm(\ref{#1})}}
\newtheorem{theorem}{\bf Theorem}[section]
\newtheorem{lemma}[theorem]{\bf Lemma}
\newtheorem{corollary}[theorem]{\bf Corollary}
\newtheorem{proposition}[theorem]{\bf Proposition}
\newtheorem{definition}[theorem]{\bf Definition}
\def\remark{\vskip 4pt \noindent {\it Remark.}\ }

\renewcommand{\theequation}{\arabic{section}.\arabic{equation}}

\begin{document}


\def\l{\lambda}
\def\p{\psi}   \def\r{\rho}
\def\s{\sigma}
\def\ve{\varepsilon}
\def\vt{\vartheta}
\def\vp{\varphi}

\def\cA{{\cal A}}
\def\bH{{\bf H}}
\def\cH{{\cal H}}
\def\cR{{\cal R}}
\def\cS{{\cal S}}

\def\el2{\ell^{\,2}}

\def\qqq{\qquad}
\def\qq{\quad}
\let\le\leqslant \let\ge\geqslant

\def\R{{\mathbb R}} \def\C{{\mathbb C}}
\def\D{{\mathbb D}} \def\T{{\mathbb T}}
\def\N{{\mathbb N}}

\def\J{\mathbb{J}}
\def\Z{\mathbb{Z}}
\def\R{\mathbb{R}}
\def\C{\mathbb{C}}
\def\T{\mathbb{T}}
\def\N{\mathbb{N}}
\def\dS{\mathbb{S}}
\def\H{\mathbb{H}}
\def\dD{\mathbb{D}}

\def\Re{\mathop{\rm Re}\nolimits}
\def\Im{\mathop{\rm Im}\nolimits}
\def\Ran{\mathop{\rm Ran}\nolimits}
\def\Iso{\mathop{\rm Iso}\nolimits}

\newcommand{\oo}[1]{{\mathop{#1}\limits^{\,\circ}}\vphantom{#1}}
\newcommand{\po}[1]{{\mathop{#1}\limits^{\phantom{\circ}}}\vphantom{#1}}
\newcommand{\nh}[1]{{\mathop{#1}\limits^{{}_{\,\bf{\wedge}}}}\vphantom{#1}}
\newcommand{\nc}[1]{{\mathop{#1}\limits^{{}_{\,\bf{\vee}}}}\vphantom{#1}}
\newcommand{\nt}[1]{{\mathop{#1}\limits^{{}_{\,\bf{\sim}}}}\vphantom{#1}}
\newcommand{\dddot}[1]{{\mathop{#1}\limits^{\bf\dots}}\vphantom{#1}}


\title {The inverse problem for perturbed harmonic oscillator on the half-line with Dirichlet boundary conditions}
\author{
Dmitry Chelkak\begin{footnote}{Dept. of Math. Analysis, Math. Mech. Faculty, St.Petersburg
State University. Universitetskij pr. 28, Staryj Petergof, 198504 St.Petersburg, Russia.
E-mail: delta4@math.spbu.ru}
\end{footnote}\ \ and\ \ Evgeny Korotyaev\begin{footnote}{
Institut f\"ur  Mathematik,  Humboldt Universit\"at zu Berlin. Rudower Chaussee 25, 12489
Berlin, Germany. 
E-mail: evgeny@math.hu-berlin.de}\end{footnote}} \maketitle

\begin{abstract}
\noindent We consider the perturbed harmonic oscillator $T_D\psi=-\psi''+x^2\psi+q(x)\psi$,
$\psi(0)=0$, in $L^2(\R_+)$, where $q\in\bH_+=\{q', xq\in L^2(\R_+)\}$ is a real-valued
potential. We prove that the mapping $q\mapsto{\rm spectral\ data}={\rm \{ eigenvalues\ of\
}T_D{\rm \}}\oplus{\rm \{norming\ constants\}}$ is one-to-one and onto. The complete
characterization of the set of spectral data which corresponds to $q\in\bH_+$ is given.
\end{abstract}

\section{Introduction and main results}

Consider the Schr\"odinger operator
\begin{equation}
\label{R3HDef} {\rm H}=-\frac{\partial^2}{\partial {\bf x}^2}+|{\bf x}|^2+q(|{\bf x}|),\ \ \
{\bf x}\in \R^3,
\end{equation}
acting in the space $L^2(\R^3)$. Let $x\!=\!|{\bf x}|$ and $q$ be a real-valued bounded
function. The operator ${\rm H}$ has pure point spectrum. Using the standard transformation
$u({\bf x})\mapsto xu(x)$ and expansion in spherical harmonics, we obtain that ${\rm H}$ is
unitary equivalent to a direct sum of the Schr\"odinger operators acting on $L^2(\R_+)$.
The first operator from this sum is given by
\begin{equation}
\label{TDDef} T_D\psi=-\psi''+x^2\psi+q(x)\psi\,,\qquad \psi(0)\!=\!0\,,\qquad x\!\ge\!0\,.
\end{equation}
The second is $-\frac{d^2}{dx^2}+x^2+\frac{2}{x^2}+q(x)$ etc. Below we consider the
simplest case, i.e., the operator $T_D$\,. In our paper we assume that
$$
q\in\bH_+=\biggl\{q\!\in\!L^2(\R_+):q',xq\!\in\!L^2(\R_+)\biggr\}\,.
$$
The similar class of potentials was used to solve the corresponding inverse problem on the real line \cite{CKK}.
Define the unperturbed operator $T_D^0\psi\!=\!-\psi''\!+\!x^2\psi$\,, $\psi(0)\!=\!0$\,. The
spectrum $\s(T_D)$ of $T_D$ is the increasing sequence of simple eigenvalues
$\sigma_{n}\!=\!\sigma_{n}^0\!+\!o(1)$\,, where $\sigma_{n}^0\!=\!4n\!+\!3$, $n\!\ge\!0$\,,
are the eigenvalues of $T_D^0$\,. Note that $\sigma(T_D)$ does not determine $q$ uniquely, see
Theorem \nolinebreak \ref{*H+DChar}. Then what does the isospectral set
$$
\Iso_D(q)=\{p\in {\bf H}_+:\, \sigma_n(p)=\sigma_n(q) {\rm\ for\ all\ } n\ge 0\}
$$
of all potentials $p$ with the same Dirichlet spectrum as $q$ look like?

The inverse problem consists of two parts:

\noindent 1) to characterize the set of all sequences of real numbers which arise as the Dirichlet spectra of $q\!\in\!\bH_+$.

\noindent 2) to describe the set $\Iso_D(q)$\,. 

We shall give the complete solution of these problems in Theorem \ref{*H+DChar}.
To  describe the set $\Iso_D(q)$ we define the {\bf norming constants}
$\nu_{n}(q)$ by\begin{footnote}{Here and below we use the notations
$\|\cdot\|_+\!=\!\|\cdot\|_{L^2(\R_+)}$\,,
$\langle\cdot,\cdot\rangle_+\!=\!\langle\cdot,\cdot\rangle_{L^2(\R_+)}$\,. } \end{footnote}
\begin{equation}
\label{NuDef} \nu_{n}(q)= \log\|\varphi(\cdot,\sigma_n(q),q)\|_{+}^{-2}=
2\log|\psi'_{n,D}(0,q)|\,,\quad n\ge 0\,,
\end{equation} 
where  $\varphi(x)=\vp(x,\l,q)$ is the solution of the
equation
\begin{equation}
\label{PertEq} -\vp''+x^2\vp+q(x)\vp=\l\vp\,, \ \
\vp(0)=0,\ \ \varphi'(0)=1,\ \ 
\quad (\l,q)\in\!\C\times \bH_+,
\end{equation}
and $\psi_{n,D}$ is the $n$-th normalized (in $L^2(\R_+)$) eigenfunction of the operator
$T_D$\,. 

\remark Let $\Psi_{n,0}({\bf x})\!=\!\Psi_{n,0}(|{\bf x}|)$ be the $n$-th normalized (in
$L^2(\R^3)$) spherically-symmetric eigenfunction of the operator $\rm H$ given by
(\ref{R3HDef}). Then we obtain
$$
\Psi_{n,0}({\bf x})=\frac{\psi_{n,D}(|{\bf x}|)}{2\sqrt{\pi}|\,{\bf x}|}\qquad {\rm and}\qquad
|\Psi_{n,0}(0)|^2=\frac{e^{\nu_{n}(q)}}{4\pi}\,,\quad n\ge 0\,.
$$

We describe papers about the inverse problem for the
perturbed harmonic oscillator, which are relevant to our paper.
McKean and Trubowitz \cite{MT}
considered the problem of reconstruction on the real line. They gave an  algorithm for the reconstruction of $q$ from norming constants for the class of real
infinitely differentiable potentials, vanishing rapidly at
$\pm\infty$\,, for fixed eigenvalues $\l_n(q)\!=\!\l_n^0$ for all $n$ and
 "norming constants" $\to\!0$ rapidly as $n\!\to\!\infty$. Later on, Levitan \cite{L} reproved some results of \cite{MT} without an exact
definition of the class of potentials. Some uniqueness theorems were obtained by Gesztesy, Simon \cite{GS1} and Chelkak, Kargaev, Korotyaev \cite{CKK1}. For uniqueness theorems
we need some asymptotics of fundamental solutions and eigenvalues
at high energy. For characterization we need "sharp"
asymptotics of these values. Usually it  is not simple.
Note that recently the asymptotics $\l_n(q)$ were determined for
bounded potentials in \cite{KKP}.
 Gesztesy and Simon  \cite{GS} proved that the each 
 $\Iso_D(q)$ is connected for various classes of
potentials. Note that the inverse problem for harmonic oscillator on the half-line for the boundary conditions $\psi'(0)=b\psi(0)$, $b\in\R$ is solved in \cite{CK}.

Our approach is based on the methods from  \cite{CKK} and \cite{PT} (devoted to the inverse Dirichlet problem on $[0,1]$).  The main point in the inverse problem for the perturbed harmonic oscillator $T_D$ on $\R_+$ is the characterization of $\Iso_D(q)$. Note that, in contrast to the case of perturbed harmonic oscillator on the real-line \cite{CKK}, the characterization of $\{\nu_n\}_{n=0}^\infty$\,, i.e. the {
parameterization of isospectral manifolds}, is given in terms of the
standard weighted $\el2$-\,space. Thus there is a big difference between
the case of the real line and the case of the half-line.

The present paper continues the series of papers \cite{CKK1}, \cite{CKK} devoted to the
inverse spectral problem for the perturbed harmonic oscillator on the real line. Note that the set of spectra which correspond to potentials from $\bH_+$ (see Sect. 3 for details) is similar to the space of spectral data in \cite{CKK}\,. In particular, the range of the linear
operator $f(z)\!\mapsto\!(1\!-\!z)^{-1/2}f(z)$ acting in some Hardy-Sobolev space in the unit disc plays an important role. As a byproduct of our analysis, we give the simple proof of the
equivalence between two definitions of the space of spectral data
(Theorem 4.2 in \cite{CKK}), which was established in
\cite{CKK} using a more complicated techniques.



We recall some basic results from \cite{CKK}. Consider the operator
$$
T\psi=-\psi''+x^2\psi+q(x)\psi\,,\quad\quad q\in\bH_{even}\!=
\!\biggl\{q\!\in\!L^2(\R):q',xq\!\in\!L^2(\R);\ q(x)\!=\!q(-x),\ x\!\in\!\R \biggr\}\,,
$$
acting in the space $L^2(\R)$. The spectrum $\s(T)$ is an increasing sequence of  simple
eigenvalues given by
$$
\l_n(q)\!=\!\l_n^0\!+\!\mu_n(q)\,,\quad {\rm where}\quad \l_n^0\!=\!\l_n(0)\!=\!2n\!+\!1,\
n\!\ge\!0\,,\quad {\rm  and}\quad  \mu_n(q)\!\to\!0\ {\rm  as} \ n\!\to\!\infty\,.
$$
Define the real weighted $\el2$-space
$$
\el2_r=\biggl\{c=\{c_n\}_{n=0}^\infty: \ c_n\!\in\! \R\,,\
\|c\|_{\el2_r}^2={\textstyle\sum_{n\ge 0}} (1\!+\!n)^{2r}|c_n|^2\!<\!+\infty\biggr\}\,, \ \
r\!\ge\!0\,,
$$
and the Hardy-Sobolev space of analytic functions in the unit disc $\D\!=\!\{z:|z|\!<\!1\}$:
$$
H^2_r=H^2_r(\D)= \biggl\{f(z)\!\equiv\!{\textstyle\sum_{n\ge 0}}f_nz^n,\ z\!\in\!\D:\
f_n\!\in\!\R\,,\ \|f\|_{H^2_r}=\|\{f_n\}_{n=0}^\infty\|_{\el2_r}\!<\!+\infty\biggr\}\,, \ \
r\!\ge\!0\,.
$$
Introduce the {\bf space of spectral data} from \cite{CKK}
\begin{equation}
\label{cHDef} \cH=\biggl\{h=\{h_n\}_{n=0}^{\infty}: \ \sum_{n\ge
0}h_nz^n\equiv\frac{f(z)}{\sqrt{1\!-\!z}}\,, \ f\!\in\!H^2_{\frac{3}{4}}\,\biggr\}\,,\quad
\|h\|_{\cH}=\|f\|_{H^2_{{3}/{4}}}\,.
\end{equation}
\begin{theorem}[\cite{CKK}]  \label{*HevenChar}
The mapping \ $q\to\{\lambda_n(q)-\lambda_n^0\}_{n=0}^\infty$ is a real-analytic isomorphism\begin{footnote}{By definition, the mapping of Hilbert spaces  $F:H_1\to H_2$ is a
local real-analytic isomorphism iff for any $y\!\in\!H_1$ it has an analytic continuation
$\widetilde{F}$ into some complex neighborhood $y\!\in\!U\!\subset\!{H_1}_\C$ of $y$ such that
$\widetilde{F}$ is a bijection between $U$ and some complex neighborhood
$F(y)\!\in\!\widetilde{F}(U)\!\subset\!{H_2}_\C$ of $F(y)$  and both
$\widetilde{F}$\,,\,$\widetilde{F}^{-1}$ are analytic. The local isomorphism $F$ is a (global)
isomorphism iff it is a bijection.}\end{footnote} between the space of even potentials
$\bH_{even}$ and the following open convex subset
$$
\cS=\biggl\{\{h_n\}_{n=0}^\infty\in\cH:
\l_0^0\!+\!h_0\!<\!\l_1^0\!+\!h_1\!<\!\l_2^0\!+\!h_2\!<\!\dots \biggr\}
\subset\cH\,.
$$
\remark {\rm The inequalities in the definition of $\cS$ correspond to the monotonicity of
eigenvalues.}
\end{theorem}

Recall  the following sharp representation from \cite{CKK}:
$$
\l_n(q)=\l_n^0+
\frac{\int_\R q(t)dt}{\pi \sqrt{\l_n^0}} + \widetilde{\mu}_n(q),\qquad
\{\widetilde{\mu}_n(q)\}_{0}^\infty\in \cH_0,\ \ q\in ????
$$
where the subspace $\cH_0\!\subset\!\cH$ of codimension $1$ is given by
\begin{equation}
\label{cH0Def} \cH_0= \biggl\{h\!\in\!\cH:\sqrt{1\!-\!z}\,
\left.{\textstyle\sum_{n\ge0}}h_nz^n\right|_{z=1}\!=\!f(1)\!=\!0\biggr\}\,.
\end{equation}
Remark that Lemma \ref{*cH0Asymp} yields
$\widetilde{\mu}_n(q)\!=\!O(n^{-\frac{3}{4}}\log^{\frac{1}{2}}n)$ as $n\!\to\!\infty$\,. Also,
we need
\begin{proposition}
\label{*q0=} {\bf (Trace formula)} For each $q\in\bH_{even}$ the following identity holds:
\begin{equation}
\label{q(0)=} q(0)=2\sum_{n\ge 0}\, (\lambda_{2n}(q)-\lambda_{2n+1}(q)+2),
\end{equation}
where the sum converges absolutely.
\end{proposition}

We come to the inverse problem for the operator $T_D$ on $\R_+$\,. For each $q\!\in\!\bH_+$ we
set $q(-x)\!=\!q(x)$\,, $x\!\ge\!0$\,. This gives a natural isomorphism between $\bH_+$ and
$\bH_{even}$\,. Then
$$
\sigma_n(q)=\lambda_{2n+1}(q),\quad n\ge 0\,.
$$
Let
\begin{equation}
\label{SDDef} \cS_D=\biggl\{\{h_n\}_{n=0}^\infty\in\cH:
\sigma_0^0\!+\!h_0\!<\!\sigma_1^0\!+\!h_1\!<\!\sigma_2^0\!+\!h_2\!<\!\dots \biggr\}\,.
\end{equation}
We formulate our main result.
\begin{theorem} \label{*H+DChar}
(i) The sequence $\{\sigma_n(q)-\sigma_n^0\}_{n=0}^\infty$ belongs to $\cS_D$ for each
potential $q\in\bH_+$\,.

\noindent (ii) For each $q\!\in\!\bH_+$ the  sequence 
$\{r_{n}(q)\}_{n=0}^\infty\in\el2_{\frac{3}{4}}$, where $r_n$ is given by
$$
\nu_{n}(q)= \nu_n^0 - \frac{q(0)}{2(2n\!+\!1)}+r_n(q),\quad and 
\quad 
\nu_n^0=\nu_n(0)=\log \biggl[\frac{4(2n\!+\!1)!}{\sqrt{\pi}\,2^{2n}[n!]^2}\biggr]\,.
$$

\noindent (iii) The mapping $q\mapsto
\left(\{\sigma_n(q)\!-\!\sigma_n^0\}_{n=0}^\infty\,,q(0)\,, \{r_{n}(q)\}_{n=0}^\infty\right)$
is a real-analytic isomorphism between $\bH_+$ and $\cS_D\times\R\times\el2_{\frac{3}{4}}$.
\remark {\rm In particular, $\left(p(0)\,,
\{r_{n}(p)\}_{n=0}^\infty\right)\in\R\times\el2_{\frac{3}{4}}$ are "independent coordinates"\
in $\Iso_D(q)$.}
\end{theorem}

The ingredients of the proof of Theorem \ref{*H+DChar} are:

\noindent i) Uniqueness Theorem. We adopt the proof from \cite{PT} and \cite{CKK1}. This proof
requires only some estimates of the fundamental solutions.

\noindent ii) Analysis of the Fr\'echet derivative of the nonlinear spectral mapping $\{{\rm
potentials}\}\!\mapsto\!\{{\rm spectral\ data}\}$ at the point $q\!=\!0$\,. We emphasize that
this linear operator is complicated (in particular, it is not the Fourier transform, as it was
in \cite{PT}). Here we essentially use the technique of generating functions (from
\cite{CKK}), which are analytic in the unit disc.

\noindent iii) Asymptotic analysis of the difference between spectral data and its Fr\'echet
derivatives at $q\!=\!0$\,. Here the calculations and asymptotics from \cite{CKK} play an
important role.

\noindent iv) The proof that the spectral mapping is a surjection, i.e. the fact that each
element of an appropriate Hilbert space can be obtained as spectral data of some potential
$q\!\in\!\bH_+$\,. Here we use the standard Darboux transform of second-order differential
equations.

\noindent {\bf The plan of the paper.}\ Sect. \ref{SectAsympt} is devoted to the basic
asymptotics of the eigenvalues $\sigma_n(q)$ and the values
$\log[(-1)^n\psi'_+(0,\sigma_n(q),q)]$\,. In Sect. \ref{SectLeadTerms} we introduce the space
$\cH$ and obtain its equivalent definition (Corollary \ref{*H0EquivDef}). Furthermore, we
consider a kind of linear approximation of our spectral data and prove Theorem \ref{*TildeQ}
that is, in a sense, the linear analogue of the main Theorem \ref{*H+DChar}. Sect.
\ref{SectAsymptNu} is devoted to the asymptotics of the norming constants $\nu_n(q)$. Also, in
this Sect. we prove Proposition \ref{*q0=}. In Sect. \ref{SectInvD} we prove the main Theorem
\ref{*H+DChar}. All needed properties of fundamental solutions, gradients of spectral data and
some technical Lemmas are collected in Appendix.




\section{\label{SectAsympt}Basic asymptotics}
\setcounter{equation}{0}

Let $\psi_+^0(x,\l)\!=\!D_{\frac{\l-1}{2}}(\sqrt{2}x)$ be the decreasing near $+\infty$
solution of the unperturbed equation
$$
-\psi''+x^2\psi = \lambda \psi
$$
We use the standard notation $D_\mu(x)$ for the Weber functions (or the parabolic cylinder
functions), see \cite{B}. Note that for each $q\!\in\!\bH_+$ the perturbed equation
$$
-\psi''+x^2\psi+q(x)\psi=\lambda\psi
$$
has the unique solution $\psi_+(x,\l,q)$ such that $\psi_+(x)\!=\!\psi_+^0(x)(1\!+\!o(1))$ as
$x\!\to\!+\infty$ (see (\ref{Psi+Asympt1})).

\begin{lemma}
\label{*NuIdent} For each $q\!\in\!\bH_+$ and $n\!\ge\!0$ the following identities hold:
\begin{equation}
\label{NuIdent1} \nu_n(q)=2\log
\biggl[-\frac{\psi'_+(0,\sigma_n(q),q)}{\dot{\psi}_+(0,\sigma_n(q),q)}\biggr],
\end{equation}
\begin{equation}
\label{NuIdent2} \dot{\psi}_+(0,\sigma_n(q),q)=
\frac{\psi_+^0(0,\sigma_n(q))}{\sigma_n(q)\!-\!\sigma_n^0}\,\,\cdot \prod_{m:m\ne
n}\frac{\sigma_n(q)\!-\!\sigma_m(q)}{\sigma_n(q)\!-\!\sigma_m^0},
\ \ \ \  \  {\partial \over\partial \lambda}{\psi}_+=\dot {\psi}_+.
\end{equation}
Remark.\ {\rm It is important that the values $\dot{\psi}_+(0,\sigma_n(q),q)$ are uniquely
determined by the spectrum $\sigma(T_D)$\,. In particular,
$\dot{\psi}_+(0,\sigma_n(p),p)\!=\!\dot{\psi}_+(0,\sigma_n(q),q)$ for all $p\!\in\!\Iso_D(q)$
and $n\!\ge\!0$\,.}
\end{lemma}

\begin{proof}
The standard identity\begin{footnote}{Here and below we use the notations $\{f,g\}=fg'-f'g$\,,
$u'=\frac{\partial}{\partial x}u$\,, $\dot u=\frac{\partial }{\partial \lambda}u$\,. }\end{footnote}
$\psi_+^2\!=\!\{\dot{\psi}_+\,,\psi_+\}'$ yields
$$
\int_0^{+\infty}\!\!\psi_+^2(x)dx= \{\dot{\psi}_+\,,\psi_+\}\Big|_{0}^{+\infty}=
-\dot{\psi}_+(0)\psi'_+(0)\,,
$$
where we omit $\sigma_n(q)$ and $q$ for short. Therefore,
$$
e^{\nu_n(q)}=[\psi_{n,D}'(0,q)]^2=
\biggl[\frac{\psi'_+(0,\sigma_n(q),q)}{\|\psi_+(\cdot,\sigma_n(q),q)\|_+}\biggr]^2=
-\frac{\psi'_+(0,\sigma_n(q),q)}{\dot{\psi}_+(0,\sigma_n(q),q)}\,.
$$
Using the Hadamard Factorization Theorem, we obtain
$$
\psi_+(0,\sigma,q)=\psi_+^0(0,\sigma)\cdot \prod_{m\ge 0}
\frac{\sigma\!-\!\sigma_m(q)}{\sigma\!-\!\sigma_m^0}\,,\qquad \sigma\in\C\,.
$$
The differentiation of $\psi_+(0,\sigma,q)$ gives (\ref{NuIdent2}).
\end{proof}

Let $\p_n^0$ be the normalized (in $L^2(\R)$) eigenfunctions of the unperturbed harmonic
oscillator on $\R$. Note that
$\psi_{n,D}^0(\cdot)\!=\!\psi_{n,D}(\cdot,0)\!=\!\sqrt{2}\psi_{2n+1}^0(\cdot)$\,. It is
well-known that
$$
\psi_n^0(x)=(n!\sqrt{\pi})^{-\frac{1}{2}}D_n(\sqrt{2}x)=
({2^n}n!\sqrt{\pi})^{-\frac{1}{2}}H_n(x)e^{-\frac{x^2}{2}}\,,\quad n\!\ge\!0\,,
$$
where $H_n(x)$ are the Hermite polynomials. For each $n\!\ge\!0$ we consider the second
solution
$$
\chi_n^0(x)=\biggl(\frac{n!\sqrt{\pi}}{2}\biggr)^{\!1/2}
\cases{(-1)^\frac{n}{2}\Im{D_{-n-1}}(i\sqrt{2}x), & $n$ is even, \cr
(-1)^\frac{n-1}{2}\Re{D_{-n-1}}(i\sqrt{2}x), & $n$ is odd,}
$$
of the equation $-\psi''\!+\!x^2\psi\!=\!\lambda_n^0 \psi$ which is uniquely defined by the
conditions
$$
\{\chi_n^0\,,\psi_n^0\}\!=\!1\,,\qquad (\psi_n^0\chi_n^0)(-x)\!=\!-(\psi_n^0\chi_n^0)(x)\,,\ \
\ x\in \R\,.
$$
Note that $(\psi_n^0\chi_n^0)(x)=(-1)^{n+1}x+O(x^2)$ as $x\!\to\!0$\,,  and $(\psi_n^0\chi_n^0)(x)=-x^{-1}+O(x^{-2})$ as $x\!\to\!\infty$,
see \cite{CKK}. Following \cite{CKK}, for $q\!\in\!\bH_+$,  we introduce
\begin{equation}
\label{NhNcQDef} \nh q_n= \langle q,(\psi_n^0)^2 \rangle_{+}\,,\qquad \nc q_n = \langle
q,\psi_n^0\chi_n^0 \rangle_{+}\,\ \ n\!\ge\!0.
\end{equation}
Also, we introduce the constants
$$
\kappa_n=\p_{+}^0(0,\l_n^0)\,,\qquad \kappa'_n=(\p_{+}^0)'(0,\l_n^0)\,,\qquad
\dot{\kappa}_n=\dot{\p}_{+}^0(0,\l_n^0)\quad {\rm and\ so\ on.}
$$

\begin{theorem}
\label{*BasicAsympt} For each $q\!\in\!\bH_+$ the
following asymptotics\begin{footnote}{Here and below $a_n\!=\!b_n\!+\!\el2_r(n)$ means that
$\{a_n\!-\!b_n\}_{n=0}^\infty\!\in\!\el2_r$\,. We say that $a_n(q)\!=\!b_n(q)\!+\!\el2_r(n)$
holds true uniformly on some set iff norms $\|\{a_n(q)\!-\!b_n(q)\}_{n=0}^\infty\|_{\el2_r}$
are uniformly bounded on this set.}\end{footnote} hold:
\begin{equation}
\label{SAsympt} \sigma_{n}(q)=\sigma_n^0+2{\nh q_{2n+1}}+\el2_{\frac{3}{4}+\delta}(n)\,,
\end{equation}
\begin{equation}
\label{OldNuAsympt} \log\frac{\psi'_+(0,\sigma_n(q),q)}{\kappa'_{2n+1}}=
\frac{\dot{\kappa}'_{2n+1}}{\kappa'_{2n+1}}\,(\sigma_n(q)\!-\!\sigma_n^0) - {\nc q_{2n+1}}+
\el2_{\frac{3}{4}+\delta}(n)\,,
\end{equation}
uniformly on bounded subsets of $\bH_+$\,, 
for some absolute constant $\delta\!>\!0$.
\end{theorem}

\noindent {\it Remark.}\ i) Proposition \ref{*NcNhasFG} immediately yields $\{\nh
q_{n}\}_{n=0}^{\infty}\!\in\!\cH$\,. Using basic properties of the spaces $\cH$, $\cH_0$ (see
Proposition \ref{*H0Emb} and Lemma \ref{*cH0Asymp}), we obtain
$$
\textstyle \nh q_{2n+1}=\pi^{-1}\int_{\R_+}q(t)dt\cdot (\sigma_n^0)^{-\frac{1}{2}}
+O(n^{-\frac{3}{4}}\log^{\frac{1}{2}}n)\,.
$$

\noindent ii) In the proof we use some technical results from \cite{CKK},
 formulated in Appendix A.1-A.4.

\begin{proof}
Let $\mu\!=\!\sigma_n(q)\!-\!\sigma_n^0$ and $m\!=\!2n\!+\!1$\,. Recall that
$\sigma_n^0\!=\!\lambda_m^0$ and $\psi_+(0,\lambda_{m}^0\!+\!\mu,q)\!=\!0$\,. Lemma
\nolinebreak \ref{*AAnalyticity} \nolinebreak (i) yields $\mu\!=\!O(m^{-1/2})$. Due to
Corollary \ref{*ADotEstim} and asymptotics (\ref{KappaValues}), we have
$$
0=\frac{\psi_+(0,\lambda_{m}^0\!+\!\mu,q)}{\dot{\kappa}_{m}} =
\frac{\psi_+^{(1)}(0,\lambda_{m}^0,q)+\dot{\kappa}_{m}\cdot\mu} {\dot{\kappa}_{m}} +
O(n^{-1}\log^2n).
$$
Hence, Lemma \ref{*LA511} (ii) gives
$
\mu = 2\nh q_{m} + O(m^{-1}\log^2m)\,.
$
Using the similar arguments, we deduce that
$$
0=\!\frac{1}{\dot{\kappa}_{m}}\,\biggl(\psi_+^{(1)} \!+ \dot{\psi}_+^0\cdot\mu + \psi_+^{(2)}
\!+ \dot{\psi}_+^{(1)}\cdot 2\nh q_{m} + \frac{\ddot{\psi}_+^0}{2}\,(2\nh q_{m})^2\biggr)
(0,\lambda_{m}^0,q) + O(m^{-\frac{3}{2}}\log^3m).
$$
Together with Lemmas \ref{*LA511} (ii), \ref{*LKappaSpecAsymp}, this yields
$$
\frac{\mu}{2}= \nh q_{m} - \nc q_{m}\nh q_{m} +
\biggl(\frac{\dot{\kappa}'_{m}}{\kappa'_{m}}\,\nh q_{m} \!+\! \frac{1}{2}\,\nc
q_{m}\biggr)\cdot 2\nh q_{m} -\frac{\ddot{\kappa}_{m}}{\kappa_{m}}\,(\nh q_{m})^2 +
\el2_{\frac{3}{4}+\delta}(m)
$$
$$
= \nh q_{m} + \biggl(\frac{2\dot{\kappa}'_{m}}{\kappa'_{m}} -
\frac{\ddot{\kappa}_{m}}{\kappa_{m}}\biggr)\cdot (\nh q_{m})^2 +\el2_{\frac{3}{4}+\delta}(m) =
\nh q_{m} + \el2_{\frac{3}{4}+\delta}(m).
$$
Furthermore, using Corollary \ref{*ADotEstim} and Lemma \ref{*LA511} (ii), we obtain
$$
\frac{\psi'_+(0,\sigma_n(q),q)}{\kappa'_{m}} = \frac{\left(\psi_+^0 + \psi_+^{(1)} +
\psi_+^{(2)}\right)(0,\lambda_{m}^0\!+\!\mu,q)} {\kappa'_{m}}+ O(m^{-\frac{3}{2}})
$$
$$
= \frac{\left( \psi_+^0 + \psi_+^{(1)} + \dot{\psi}_+^0\cdot\mu + \psi_+^{(2)} +
\dot\psi_+^{(1)}\cdot 2\nh q_{m} + \frac{1}{2}\,\ddot\psi_+^0\cdot(2\nh q_{m})^2
\right)'\!(0,\lambda_{m}^0,q)}{\kappa'_{m}}+ O(m^{-\frac{3}{2}}\log^3m)
$$
$$
= 1- \nc q_{m} + \frac{\dot{\kappa}'_{m}}{\kappa'_{m}}\,\cdot \mu + \biggl( \frac{1}{2}\,(\nc
q_{m})^2 \!-\! \frac{\pi^2}{8}\,(\nh q_{m})^2 \biggr) + \biggl(-
\frac{\dot{\kappa}'_{m}}{\kappa'_{m}}\,\nc q_{m} \!+\! \frac{\pi^2}{8}\,\nh q_{m}\biggr)\cdot
2\nh q_{m} + \frac{2\ddot{\kappa}'_{m}}{\kappa'_{m}}\,(\nh q_{m})^2+
\el2_{\frac{3}{4}+\delta}(m).
$$
Hence,
$$
\log\frac{\psi'_+(0,\sigma_n(q),q)}{\kappa'_{2n+1}} = -\nc q_{m} +
\frac{\dot{\kappa}'_{m}}{\kappa'_{m}}\,\cdot \mu
+2\biggl(\frac{\ddot{\kappa}'_{m}}{\kappa'_{m}}-
\frac{(\dot{\kappa}'_{m})^2}{(\kappa'_{m})^2}+\frac{\pi^2}{16}\biggr)\cdot (\nh
q_{m})^2+\el2_{\frac{3}{4}+\delta}(m)
$$
$$
=-\nc q_{m} + \frac{\dot{\kappa}'_{m}}{\kappa'_{m}}\,\cdot \mu +\el2_{\frac{3}{4}+\delta}(m),
$$
where we have used Lemma \ref{*LKappaSpecAsymp}.
\end{proof}

\section{\label{SectLeadTerms}Coefficients $\bf {\mathop{q}\limits^\wedge}_{2n+1}$\,,
$\bf {\mathop{q}\limits^\vee}_{2n+1}$ and ${\bf \mathop{q}\limits^\sim}_n$}
\setcounter{equation}{0}

Let
$$
\widetilde{\psi}_n^0(x)\!=\! 2^{1/4}\psi_n^0(\sqrt{2}x)\,,\quad n\!\ge\!0\,.
$$
Note that the mapping
\begin{equation}
\label{x31} q\mapsto \{\langle q,\widetilde{\psi}_n^0 \rangle \}_{n=0}^\infty\,,\qquad \bH\to
\el2_{1/2}\,.
\end{equation}
is a linear isomorphism\begin{footnote}{We say that the linear operator is a linear
isomorphism iff it is bounded and its inverse is bounded too.}\end{footnote}. Moreover, since
$\{\widetilde{\psi}_{2m}^0\}_{m=0}^\infty$ is the orthogonal basis of the space
$\bH_{even}$\,, it is the orthogonal basis of $\bH_+$\,. On the contrary,
$\{\widetilde{\psi}_{2m+1}^0\}_{m=0}^\infty$ is the orthogonal basis of the subspace
$$
\bH_+^0=\{q\!\in\!\bH_+:q(0)\!=\!0\}\subsetneq\bH_+\,.
$$

Following \cite{CKK}, for each potential $q\in \bH_+$ we define two (analytic in the unit disc
$\D$) functions
\begin{equation}\label{FqGqDef}
\begin{array}{c}\displaystyle
(Fq)(z)\equiv\frac{1}{(2\pi)^{1/4}}\sum_{k\ge 0} \sqrt{E_k}\,\langle q,\widetilde{\psi}_{2k}^0
\rangle_+ \cdot z^k\,,\quad z\!\in\!\D\,,\cr \displaystyle
(Gq)(z)\equiv-\frac{(2\pi)^{1/4}}{2}\sum_{k\ge 0} \frac{\langle
q,\widetilde{\psi}_{2k+1}^0\rangle_+}{\sqrt{(2k\!+\!1)E_k}}\,z^k\,,\quad z\!\in\!\D\,,
\end{array}
\end{equation}
where
\begin{equation}
\label{EAsympt} E_k=\frac{(2k)!}{2^{2k}(k!)^2}\sim \pi^{-\frac{1}{2}}k^{-\frac{1}{2}}\quad
{\rm as}\quad k\!\to\!\infty\,.
\end{equation}

\begin{lemma} \label{L31}
(i) The mapping $q\mapsto Fq$ is a linear isomorphism between $\bH_+$ and $H^2_{3/4}$\,.

\noindent (ii) The mapping $q\mapsto Gq$ is a linear isomorphism between $\bH_+^0$ and
$H^2_{3/4}$\,.

\noindent (iii) $(Gq)(\cdot)\in C(\T\setminus\{-1\})$ and $(Gq)(\cdot)\in L^1(\T)$ for each
$q\in\bH_+$.
\end{lemma}

\begin{proof}
(i), (ii) Using (\ref{x31}) and (\ref{EAsympt}), we deduce that the mappings
$$
q \mapsto \{\langle q, \widetilde{\psi}_{2m}^0 \rangle_+ \}_{m=0}^\infty \mapsto Fq\,,\qquad
\bH_+ \to \el2_{1/2} \to H^2_{3/4}
$$
$$
q \mapsto \{\langle q, \widetilde{\psi}_{2m+1}^0 \rangle_+ \}_{m=0}^\infty \mapsto Gq\,,\qquad
\bH_+^0 \to \el2_{1/2} \to H^2_{3/4}
$$
are linear isomorphisms.

\noindent (iii) Let $q\in \bH_+$\,. Since
$\widetilde{\psi}_0^0(0)\!=\!2^{\frac{1}{4}}\pi^{-\frac{1}{4}}$ (see (\ref{x31_1})), we obtain
$$
q(x)=2^{-\frac{1}{4}}\pi^{\frac{1}{4}}q(0)\cdot \widetilde{\psi}_0^0(x)+ q_0(x)\,,\quad
x\!\ge\! 0\,,
$$
for some $q_0\!\in\!\bH_+^0$\,. Due to (ii), we have $Gq_0\in H^2_{3/4}\subset C(\T) \subset
L^1(\T)$. Furthermore, (\ref{x31_2}) yields
$$
(G\widetilde{\psi}_0^0)(z) = -\frac{(2\pi)^{1/4}}{2\sqrt{2\pi}} \sum_{m\ge 0} \frac {(-1)^m
z^m }{2m\!+\!1}\,.
$$
Hence, $G\widetilde{\psi}_0^0\in C(\T\setminus\{-1\})$, $G\widetilde{\psi}_0^0\in L^1(\T)$ and
the same holds for $Gq$\,.
\end{proof}


\begin{lemma}[\cite{CKK}]
\label{*NcNhasFG} Let $q\!\in\!\bH_+$. Then the following identities\begin{footnote}{We write
$f(z)\equiv g(z)$ iff the identity $f(z)=g(z)$ holds true for all
$z\!\in\!\D$\,.}\end{footnote} hold:
\begin{equation}
\label{NcNhAsFqGq} \sum_{n\ge 0} \nh q_nz^n \equiv\frac{(Fq)(z)}{\sqrt{1\!-\!z}}\,,\qquad
\sum_{n\ge 0} \nc q_nz^n \equiv
P_+\biggl[\frac{(Gq)(\zeta)}{\sqrt{1\!-\!\overline{\zeta}}}\biggr]\,,\qquad z\in\D\,,
\end{equation}
\begin{equation}
\label{Fq1-1} (Fq)(1)\!=\!(2\pi)^{-\frac{1}{2}}\int_{\R_+}q(t)dt,\qquad
(Fq)(-1)\!=\!2^{-\frac{3}{2}}q(0),
\end{equation}
where the coefficients $\nh q_n$ and $\nc q_n$\,, $n\!\ge\!0$\,, are defined by (\ref{NhNcQDef}).
\remark {\rm i) Here and below we put $(P_+f)(z)\equiv\frac{1}{2\pi i}
\int_{|\zeta|=1}\frac{f(\zeta)d\zeta}{\zeta-z}$ for any $f\!\in\!L^1(\T)$ and $z\!\in\!\D$\,.
In particular, the identity $(P_+\sum_{n=-k}^k c_n{\zeta}^n)(z)\!\equiv\!\sum_{n=0}^kc_nz^n$
holds true for any $c_n\!\in\!\C$\,.

\noindent ii) Definition (\ref{cHDef}) of the space $\cH$ is directly motivated by asymptotics
(\ref{SAsympt}) and Eq. (\ref{NcNhAsFqGq}).}
\end{lemma}
\begin{proof}
Identities (\ref{NcNhAsFqGq}) were proved in \cite{CKK} (Propositions 1.2 and 2.9). Also, in
\cite{CKK} it was shown that
$$
1=\sum_{k\ge 0}\widetilde\psi_{2k}^0(x)\int_\R\!\widetilde{\psi}_{2k}^0(t)dt =
(2\pi)^{\frac{1}{4}}\sum_{k\ge 0}\sqrt{E_k}\,\widetilde{\psi}_{2k}^0(x)\,,
$$
in the sense of distributions, which gives
$(Fq)(1)\!=\!(2\pi)^{-\frac{1}{2}}\int_{\R_+}q(t)dt$\,. Furthermore,
$$
\delta(x)=\sum_{k\ge 0}\widetilde{\psi}_{2k}^0(x)\cdot\widetilde\psi_{2k}^0(0) = \sum_{k\ge
0}\widetilde{\psi}_{2k}^0(x)\cdot \frac{2^{1/4}H_{2k}(0)}{(\sqrt{\pi}\,2^{2k}(2k)!)^{1/2}}=
\biggl(\frac{2}{\pi}\biggr)^{\!1/4}\sum_{k\ge 0}(-1)^k\sqrt{E_k}\,\widetilde{\psi}_{2k}^0(x)
$$
in the sense of distributions. Together with (\ref{FqGqDef}), this implies
$(Fq)(-1)\!=\!2^{-\frac{1}{2}}\cdot\frac{q(0)}{2}$\,.
\end{proof}

We need some results from \cite{CKK} (see Lemmas 2.10, 2.11 \cite{CKK}).

\begin{proposition}
\label{*H0Emb} (i) For each $\{h_n\}_{n=0}^\infty\!\in\!\cH$ there exist unique $v\!\in\!\R$
and $\{h_n^{(0)}\}_{n=0}^\infty\!\in\!\cH_0$ such that
$h_n\!=\!{v}\cdot(\l_n^0)^{-\frac{1}{2}}+h_n^{(0)}$. The mapping
$$
h\mapsto(v,h^{(0)})
$$
is a linear isomorphism between $\cH$ and $\R\times\cH_0$\,. If $h\!=\!\{\nh
q_n\}_{n=0}^\infty$\,, $q\!\in\!\bH_+$\,, then
$$
\textstyle
v\!=\!2^{\frac{1}{2}}\pi^{-\frac{1}{2}}\sqrt{1\!-\!z}\,\sum_{n=0}^{+\infty}h_nz^n\Big|_{z=1}\!=
\!\pi^{-1}\!\int_{\R_+} q(t)dt\,.
$$

\noindent  (ii) The set of finite sequences $\left\{(h_0\,,\dots\,,h_k\,,0\,,0\,,\dots),\
k\!\ge\!0\,,\ h_j\!\in\!\R\right\}$ is dense in $\cH_0$\,.

\noindent (iii) The embeddings $\el2_{3/4}\!\subset\!\cH_0\!\subset\!\el2_{1/4}$ are
fulfilled.

\noindent (iv) If $\{h_n\}_{n=0}^\infty\!\in\!\cH$, then
$\{h_n-h_{n+1}\}_{n=0}^\infty\!\in\!\ell^2_{3/4}$.

\remark {\rm Since $\cH_0\!\subset\!\el2_{1/4}$\,, the sequence of leading terms
$\{(\l_n^0)^{-\frac{1}{2}}v\}_{n=0}^{\infty}$ doesn't belong to $\cH_0$\,.}

\end{proposition}

The next Lemma gives the $O$-type estimate for sequences from $\cH_0$\,.

\begin{lemma} \label{*cH0Asymp} Let $\{h_n\}_{n=0}^\infty\!\in\!\cH_0$\,. Then
$h_n=O(n^{-\frac{3}{4}}\log^{\frac{1}{2}}n)$ as $n\!\to\!\infty$\,.
\end{lemma}
\begin{proof}
The proof is similar to the proof of Lemma 2.1 in \cite{Ch}.
Definition (\ref{cH0Def}) of $\cH_0$ yields
$$
\sum_{n\ge 0}h_nz^n\equiv\frac{\sum_{k\ge 0}f_kz^k}{\sqrt{1\!-\!z}}\,,\qquad
\{f_k\}_{k=0}^{\infty}\!\in\!\el2_{\frac{3}{4}}\,,\qquad \sum_{k\ge 0}f_k\!=\!f(1)\!=\!0\,.
$$
Recall that $(1\!-\!z)^{-\frac{1}{2}}\!\equiv\!\sum_{m\ge 0}E_mz^m$\,. Hence,
$$
h_n=\sum_{k=0}^nE_{n-k}f_k=\sum_{k=1}^n(E_{n-k}\!-\!E_n)f_k-E_n\sum_{k=n+1}^{\infty}f_k\,.
$$
It is easy to see that $E_n\!=\!O(n^{-\frac{1}{2}})$ and
$E_{n-k}\!-\!E_n\!=\!O(kn^{-1}(n\!-\!k\!+\!1)^{-\frac{1}{2}})$\,. Therefore,
$$
\biggl|E_n\sum_{k=n+1}^{\infty}f_k\biggr| \le
E_n\biggl(\sum_{k=n+1}^{+\infty}k^{-\frac{3}{2}\,}\biggr)^{\!\!1/2}
\biggl(\sum_{k=n+1}^{+\infty}k^{\frac{3}{2}}|f_k|^2\biggr)^{\!\!1/2}=O(n^{-\frac{3}{4}})\,,
$$
$$
\biggl|\sum_{k=1}^n(E_{n-k}\!-\!E_n)f_k\biggr|\le \biggl(\sum_{k=1}^n
O\biggl(\frac{k^{\frac{1}{2}}}{n^2(n\!-\!k\!+\!1)}\biggr)\,\biggr)^{\!\!1/2}
\biggl(\sum_{k=1}^nk^{\frac{3}{2}}|f_k|^2\biggr)^{\!\!1/2} =
O(n^{-\frac{3}{4}}\log^{\frac{1}{2}}n)\,,
$$
where the estimate $\sum_{k=1}^{+\infty}k^{\frac{3}{2}}|f_k|^2\!<\!+\infty$ has been used.
\end{proof}

Recall that $Fq$ and $Gq$ are defined by (\ref{FqGqDef}) and the system of functions
$\{\widetilde{\psi}_{2k}^0\}_{k=0}^\infty$ is a basis of $\bH_+$. Therefore, it is possible to
rewrite $Gq$ in terms of $Fq$\,. Note that this situation differs from the case $q\!\in\!\bH$
(the perturbed harmonic oscillator on the whole real line \cite{CKK}), where the functions
$Fq$ and $Gq$ are "independent coordinates"\ in the space of potentials.

For $\zeta=e^{i\phi}\in\!\T, \phi\!\in\!(-\pi,\pi)$\,, $\zeta\!\ne\!-1$\,, we define
$
\sqrt{\zeta}\!=\!e^{\frac{i\phi}{2}}.
$
We have the identity
\begin{equation}
\label{SqrtDef}
 \frac{1}{\sqrt{\zeta}}=
\frac{2}{\pi}\sum_{s\in \Z}\frac{(-1)^s}{2s\!+\!1}\,\zeta^s \quad {\rm in}\ \
L^2(\T)\,.
\end{equation}
\begin{lemma}
\label{*GasF} For each $q\!\in\!\bH_+$ the following identity holds:
\begin{equation}
\label{GasF} (Gq)(z)\equiv
-\frac{\pi}{2}P_+\biggl[\frac{(Fq)(\zeta)}{\sqrt{\zeta}}\biggr]\,,\quad z\!\in\!\D\,.
\end{equation}
\end{lemma}
\begin{proof}
We determine the coefficients of the function $\p_{2m+1}^0$ with respect to the basis
$\{\p_{2k}^0\}_{k=0}^\infty$. The standard identity $\{\psi_{2k}^0\,,\psi_{2m+1}^0\}'\!=\!
(\lambda_{2k}^0\!-\!\lambda_{2m+1}^0)\psi_{2k}^0\psi_{2m+1}^0$ yields
$$
({\p}_{2m+1}^0\,,{\psi}_{2k}^0)_{+}= \int_{\R_+}\p_{2m+1}^0(x)\psi_{2k}^0(x)dx=
\frac{\{\p_{2k}^0\,,\p_{2m+1}^0\}(0)}{\l_{2m+1}^0\!-\!\l_{2k}^0}=
\frac{\p_{2k}^0(0)(\p_{2m+1}^0)'(0)}{2(2(m\!-\!k)\!+\!1)}\,.
$$
Note that
\begin{equation}
\label{x31_1} \psi_{2k}^0(0)=\frac{H_{2k}(0)}{(\sqrt{\pi}\,2^{2k}(2k)!)^{1/2}}=
\frac{(-1)^k2^k(2k\!-\!1)!!}{(\sqrt{\pi}\,2^{2k}(2k)!)^{1/2}}=
\frac{(-1)^k}{\pi^{1/4}}\,\sqrt{E_k}\,,
\end{equation}
$$
(\psi_{2m+1}^0)'(0)=\frac{H_{2m+1}'(0)}{(\sqrt{\pi}\,2^{2m+1}(2m\!+\!1)!)^{1/2}}=
\frac{(-1)^m2^{m+1}(2m\!+\!1)!!}{(\sqrt{\pi}\,2^{2m+1}(2m\!+\!1)!)^{1/2}}=
\frac{(-1)^m\sqrt{2}}{\pi^{1/4}}\,\sqrt{(2m\!+\!1)E_m}\,.
$$
Therefore,
\begin{equation}
\label{x31_2} (\widetilde{\p}_{2m+1}^0\,,\widetilde{\psi}_{2k}^0)_{+}=
({\p}_{2m+1}^0\,,{\psi}_{2k}^0)_{+}=
\frac{1}{\sqrt{2\pi}}\cdot\frac{(-1)^{m-k}}{2(m\!-\!k)\!+\!1}\,\sqrt{(2m\!+\!1)E_mE_k}\,.
\end{equation}
Since $\|\widetilde{\psi}_{2k}^0\|_{+}^2\!=\!\frac{1}{2}$\,, we obtain
$$
\frac{\widetilde{\psi}_{2m+1}^0}{\sqrt{(2m\!+\!1)E_m}}=
\sqrt{\frac{2}{\pi}}\,\sum_{k=0}^{+\infty} \frac{(-1)^{m-k}}{2(m\!-\!k)\!+\!1}\cdot
\sqrt{E_k}\,\widetilde{\psi}_{2k}^0\,.
$$
This gives
$$
P_+\biggl[\frac{(Fq)(\zeta)}{\sqrt{\zeta}}\biggr]\equiv \frac{2^{3/4}}{\pi^{5/4}}\,
P_+\biggl[\,\sum_{l=-\infty}^{+\infty}\frac{(-1)^l}{2l\!+\!1}\,\zeta^l\cdot
\sum_{k=0}^{+\infty}\sqrt{E_k}(q,\widetilde\psi_{2k}^0)_{+}\zeta^k\,\biggr]
$$
$$
\equiv \frac{2^{1/4}}{\pi^{3/4}}\,\sum_{m=0}^{+\infty}
\frac{(q,\widetilde{\psi}_{2m+1}^0)_{+}}{\sqrt{(2m\!+\!1)E_m}}\,z^m \equiv
-\frac{2}{\pi}\,(Gq)(z)\,,\quad z\!\in\!\D\,,
$$
where definition (\ref{FqGqDef}) of the functions $Fq$ and $Gq$ has been used.
\end{proof}

 We introduce the formal linear operator $\cA$ by
\begin{equation}
\label{OpSqrtZeta} (\cA f)(z)\equiv P_+\biggl[\frac{f(\zeta)}{\sqrt{-\zeta}}\biggr] \equiv
P_+\biggl[\sqrt{1\!-\!\overline{\zeta}}\cdot\frac{f(\zeta)}{\sqrt{1\!-\!\zeta}}\biggr]\,,
\qquad f\!\in\!L^1(\T)\,.
\end{equation}
Let
$$
\oo{H}^2_{3/4}=\{f\!\in\!\po{H}^2_{3/4}:f(1)\!=\!0\}\subset \po{H}^2_{3/4}\,.
$$

Using Lemma \ref{*GasF} we shall obtain the simple proof of  Theorem 4.2
from \cite{CKK} about the equivalent definition of $\cH_0$\,.

\begin{corollary}\label{*H0EquivDef} (i) The operator $\cA:\oo H^2_{3/4}\to\po H^2_{3/4}$ and
its inverse are bounded.\\
(ii) The following identity holds:
\begin{equation}
\label{cH0Def2} \cH_0=\biggl\{\{h_n\}_{n=0}^\infty: \sum_{n\ge 0}h_nz^n\equiv
P_+\biggl[\frac{g(\zeta)}{\sqrt{1\!-\!\overline{\zeta}}}\biggr],\ g\!\in\!H^2_{3/4}\biggr\}\,.
\end{equation}
The norms $\|h\|_{\cH}$ and $\|g\|_{H^2_{3/4}}$ are equivalent, i.e. $C_1\|g\|_{H^2_{3/4}}\le
\|h\|_{\cH}\le C_2\|g\|_{H^2_{3/4}}$ for any $g\in H^2_{3/4}$ and some absolute constants
$C_1\,,C_2>0$. \remark {\rm This equivalence was proved in \cite{CKK} using different and
complicated arguments.}
\end{corollary}
\begin{proof} (i) Recall that the mapping
$q\mapsto (Gq)(-z)$ is a linear isomorphism between $\bH_+^0$ and $H^2_{3/4}$\,. Also, due to
the identity $(Fq)(-1)\!=\!2^{-3/2}q(0)$ (see Lemma \ref{*NcNhasFG}), the mapping $q\mapsto
(Fq)(-z)$ is a linear isomorphism between $\bH_+^0$ and $\oo{H}^2_{3/4}$\,. Therefore, the
mapping
$$
f(z)\equiv(Fq)(-z)\mapsto q \mapsto (Gq)(-z)\equiv
-\frac{\pi}{2}\,P_+\biggl[\frac{(Fq)(-\zeta)}{\sqrt{\zeta}}\biggr]\equiv \cA f(z)
$$
is a linear isomorphism between $\oo H^2_{3/4}$\,, $\po \bH_+^0$ and $\po H^2_{3/4}$
respectively.\\
(ii) If $g\!\in\!H^2_{3/4}$\,, then $g_0(z)\!\equiv\!g(z)\!-\!g(1)\!\in\!\oo H^2_{3/4}$ and so
$|g_0(\zeta)|\!\le\!C|\zeta\!-\!1|^{1/4}$, $|\zeta|\!=\!1$\,, for some constant $C\!>\!0$\,.
Hence, the following equivalence is valid:
$$
\begin{array}{c}\displaystyle
\sum_{n\ge 0}h_nz^n\equiv P_{+}\biggl[\frac{g(\zeta)}{\sqrt{1\!-\!\overline{\zeta}}}\biggr] \
\ \Leftrightarrow\ \ \sum_{n\ge 0}h_nz^n\equiv g(1)+\frac{g_0(z)}{\sqrt{1\!-\!\overline{z}}}-
P_{-}\biggl[\frac{g_0(\zeta)}{\sqrt{1\!-\!\overline{\zeta}}}\biggr]
\cr\displaystyle\Leftrightarrow\ \ P_{+}\biggl[\sqrt{1\!-\!\overline{\zeta}}\sum_{n\ge 0}
h_n\zeta^n\biggr] {\vphantom{\Big|^\big|}}\equiv g(1)+g_0(z)\equiv g(z)\,,
\end{array}
$$
where $P_{-}f\!\equiv\!f\!-\!P_{+}f$ is the projector to the subspace of antianalytic
functions in $\D$\,. Therefore, the equation
$$
\frac{f(z)}{\sqrt{1\!-\!z}}\equiv\sum_{n\ge 0}h_nz^n\equiv
P_+\biggl[\frac{g(\zeta)}{\sqrt{1\!-\!\overline{\zeta}}}\biggr]\,,\quad {\rm where}\quad
f\!\in\!\oo{H}^2_{3/4}\,,\ \ g\!\in\!\po{H}^2_{3/4}\,,
$$
is equivalent to $g(z)\equiv(\cA f)(z)$\,. Then, (\ref{cH0Def2}) follows from (i).
\end{proof}

\begin{lemma}
\label{*NcAsNh} For each $q\!\in\!\bH_+$ the following identity holds:
\begin{equation}
\label{NcAsNh} \sum_{n\ge 0}\nc q_{2n+1}z^n \equiv \frac{\pi}{2}\,P_+\biggl[
\frac{1}{\sqrt{-\zeta}} \sum_{n\ge 0} \nh q_{2n} \zeta^n\biggr]\,,\quad z\!\in\!\D\,.
\end{equation}
\end{lemma}
\begin{proof} Due to identities (\ref{FqGqDef}) and Lemma \ref{*GasF}, we have
$$
\sum_{n\ge 0} \nc q_n z^n \equiv P_+\biggl[ \frac{(Gq)(\zeta)}{\sqrt{1\!-\!\overline{\zeta}}}
\biggr] \equiv -\frac{\pi}{2} P_+\biggl[
\frac{(Fq)(\zeta)}{\sqrt{\zeta}\sqrt{1\!-\!\overline{\zeta}}} \biggr]
$$
$$
\equiv -\frac{\pi}{2} P_+\biggl[
\frac{\sqrt{1\!-\!\zeta}}{\sqrt{\zeta}\sqrt{1\!-\!\overline{\zeta}}} \sum_{n\ge 0} \nh q_n
\zeta^n \biggr] \equiv -\frac{\pi}{2} P_+\biggl[ \frac{\sqrt{-\zeta}}{\sqrt{\zeta}} \sum_{n\ge
0} \nh q_n \zeta^n \biggr]\,.
$$
Therefore,
$$
\sum_{n\ge 0} \nc q_{2n+1} z^{2n+1} \equiv -\frac{\pi}{4} P_+\biggl[
\frac{\sqrt{-\zeta}}{\sqrt{\zeta}} \sum_{n\ge 0} \nh q_n \zeta^n -
\frac{\sqrt{\zeta}}{\sqrt{-\zeta}} \sum_{n\ge 0} \nh q_n (-\zeta)^n \biggr]
$$
$$
\equiv \frac{\pi}{4} P_+\biggl[ \frac{\sqrt{\zeta}}{\sqrt{-\zeta}}\sum_{n\ge 0} \nh q_n
\left(\zeta^n\!+\!(-\zeta)^n\right) \biggr] \equiv \frac{\pi}{2} P_+\biggl[
\frac{\sqrt{\zeta}}{\sqrt{-\zeta}}\sum_{n\ge 0} \nh q_{2n} \zeta^{2n} \biggr]\,.
$$
This yields
$$
\sum_{n\ge 0} \nc q_{2n+1} z^{2n} \equiv \frac{\pi}{2} P_+\biggl[
\frac{1}{\zeta}\cdot\frac{\sqrt{\zeta}}{\sqrt{-\zeta}}\sum_{n\ge 0} \nh q_{2n} \zeta^{2n}
\biggr] \equiv \frac{\pi}{2} P_+\biggl[ \frac{1}{\sqrt{-\zeta^2}}\sum_{n\ge 0} \nh q_{2n}
\zeta^{2n} \biggr]\,,
$$
since $\sqrt{-\zeta}\cdot\sqrt{\zeta}=\sqrt{-\zeta^2}$ for $\zeta\!\in\!\T$\,,
$\zeta\!\ne\!\pm 1$\,.
\end{proof}


We consider linear terms $\{\nh q_{2n+1}\}_{n=0}^\infty$ and $\{\nc q_{2n+1}\}_{n=0}^\infty$
in asymptotics (\ref{SAsympt}), \nolinebreak (\ref{OldNuAsympt}).

\begin{proposition}
\label{*PropNhNC} (i) For each $q\!\in\!\bH_+$ the following identity is fulfilled:
\begin{equation}
\label{NhAsFD} \sum_{n\ge 0} \nh q_{2n+1} z^n \equiv \frac{(F_Dq)(z)}{\sqrt{1\!-\!z}}\,,\quad
z\!\in\!\D\,,
\end{equation}
where
$$
(F_Dq)(z^2)\equiv
\frac{1}{2z}\left((Fq)(z)\sqrt{1\!+\!z}-(Fq)(-z)\sqrt{1\!-\!z}\right).
$$

\noindent (ii) For each $q\!\in\!\bH_+$ the following identity is fulfilled:
\begin{equation}
\label{NcAsGD} \sum_{n\ge 0} \nc q_{2n+1} z^n \equiv \frac{\pi}{2} P_+\biggl[
\frac{1}{\sqrt{-\zeta}} \biggl( (G_Dq)(\zeta) + \sum_{n\ge 0} \nh q_{2n+1} \zeta^n  \biggr)
\biggr]\,,\quad z\!\in\!\D\,,
\end{equation}
where
$$
(G_Dq)(z^2)\equiv
\frac{1}{2z}\left((Fq)(-z)\sqrt{1\!+\!z}-(Fq)(z)\sqrt{1\!-\!z}\right).
$$

\noindent (iii) The mapping
$$
q\mapsto (F_Dq\,;G_Dq)
$$
is a linear isomorphism between $\bH_+$ and $H^2_{3/4}\times H^2_{3/4}$\,.
\end{proposition}

\begin{proof} (i) Due to (\ref{FqGqDef}), we have
$$
\sum_{n\ge 0} \nh q_{2n+1} z^{2n} \equiv \frac{1}{2z}\biggl( \frac{(Fq)(z)}{\sqrt{1\!-\!z}} -
\frac{(Fq)(-z)}{\sqrt{1\!+\!z}} \biggr) \equiv \frac{(F_Dq)(z^2)}{\sqrt{1\!-\!z^2}}\,.
$$
This gives (\ref{NhAsFD}).

\noindent (ii) Recall that Lemma \ref{*NcAsNh} yields
$$
\sum_{n\ge 0}\nc q_{2n+1}z^{n} \equiv \frac{\pi}{2}\,P_+\biggl[ \frac{1}{\sqrt{-\zeta}}
\sum_{n\ge 0} \nh q_{2n} \zeta^{n}\biggr]\,.
$$
Using (\ref{FqGqDef}), we obtain
$$
\sum_{n\ge 0} (\nh q_{2n} - \nh q_{2n+1}) z^{2n} \equiv \frac{1}{2}\biggl(
\frac{(Fq)(z)}{\sqrt{1\!-\!z}} + \frac{(Fq)(-z)}{\sqrt{1\!+\!z}} \biggr) -
\frac{(F_Dq)(z^2)}{\sqrt{1\!-\!z^2}} \equiv (G_Dq)(z^2)\,.
$$
This gives (\ref{NcAsGD}).

\noindent (iii) Recall that the mapping $q\mapsto Fq$ is a linear isomorphism between $H_+$
and $H^2_{3/4}$\,. Therefore, we need to prove that $Fq\mapsto(F_Dq\,;G_Dq)$ is a linear
isomorphism between $H^2_{3/4}$ and $H^2_{3/4}\times H^2_{3/4}$\,. Due to definitions of $F_D$
and $G_D$\,, the direct mapping is bounded. Since
$$
(Fq)(z)\equiv (F_Dq)(z^2)\sqrt{1\!+\!z}+(G_Dq)(z^2)\sqrt{1\!-\!z}\,,
$$
the inverse mapping is bounded too.
\end{proof}

\begin{definition} \label{*NtDef} For $q\!\in\!\bH_+$ define coefficients $\nt q_n$\,, $n\!\ge\!0$\,, by
$$
\sum_{n\ge 0} \nt q_n z^n \equiv \frac{\pi}{2} P_+\biggl[ \frac{1}{\sqrt{-\zeta}} \biggl(
(G_Dq)(\zeta) - (G_Dq)(1) \biggr) \biggr]\,,\quad z\!\in\!\D\,.
$$
Remark.\ {\rm Due to Proposition \ref{*PropNhNC}, we have $G_Dq\!\in\!\po H^2_{3/4}$\,. Hence,
$(G_Dq)(\cdot)\!-\!(G_Dq)(1)\!\in\!\oo H^2_{3/4}$ and Corollary \ref{*H0EquivDef} gives $\{\nt
q_n\}_{n=0}^\infty\!\in\!\el2_{3/4}$\,.}
\end{definition}
\begin{theorem}
\label{*TildeQ} (i) For each $q\!\in\!\bH_+$ the following identities hold:
\begin{equation}
\label{Nc=Nt} \nc q_{2n+1}= \frac{q(0)}{4(2n\!+\!1)}+ \nt q_n + \sum_{m\ge 0} \frac{\nh
q_{2m+1}}{2(n\!-\!m)+1}\,,\quad n\!\ge\!0\,.
\end{equation}
(ii) The mapping
$$
q\mapsto \left(\{\nh q_{2n+1}\}_{n=0}^\infty\,; q(0)\,; \{\nt q_n\}_{n=0}^\infty\right)
$$
is a linear isomorphism between $\bH_+$ and $\cH\times\R\times\el2_{3/4}$\,.
\end{theorem}

\begin{proof}
(i) Due to identity (\ref{NcAsGD}) and Definition \ref{*NtDef}, we have
\begin{equation}
\label{xNc=Nt} \sum_{n\ge 0} \nc q_{2n+1}z^n \equiv \sum_{n\ge 0} \nt q_n z^n + \frac{\pi}{2}
P_+\biggl[ \frac{1}{\sqrt{-\zeta}} \biggl( (G_Dq)(1) + \sum_{n\ge 0} \nh q_{2n+1} z^n \biggr)
\biggr]\,.
\end{equation}
Note that $(G_Dq)(1)=2^{-\frac{1}{2}}(Fq)(-1)$\,. Then, identity (\ref{Fq1-1}) yields
$(G_Dq)(1)= \frac{1}{4}{q(0)}$\,.
 Substituting the identity
$
\frac{1}{\sqrt{-\zeta}}= \frac{2}{\pi}\sum_{s\in \Z}\frac{\zeta^s}{2s\!+\!1} \quad
{\rm in}\ \ L^2(\T)
$
(see (\ref{SqrtDef})) into (\ref{xNc=Nt}), we obtain (\ref{Nc=Nt}).

\noindent (ii) Due to Proposition \ref{*PropNhNC} and identity (\ref{NhAsFD}), the mappings
$$
\begin{array}{ccccc}
\vphantom{\big|_|} q & \mapsto & (F_Dq\,;G_Dq) & \mapsto & (\{\nh q_{2n+1}\}_{n=0}^\infty\,;
G_Dq)\,, \cr \vphantom{\big|^|}  \bH_+ & \to & H^2_{3/4}\times H^2_{3/4} & \to & \cH\times
H^2_{3/4}\,,
\end{array}
$$
are linear isomorphisms. Using Corollary \ref{*H0EquivDef}, we deduce that the mapping
$$
G_Dq \mapsto ((G_Dq)(1)\,;\{\nt q_n\}_{n=0}^\infty)\,,\qquad H^2_{3/4}\to \R\times
\el2_{3/4}\,,
$$
is a linear isomorphism too. The identity $(G_Dq)(1)= \frac{1}{4}{q(0)}$ completes the proof.
\end{proof}

\section{\label{SectAsymptNu}Asymptotics of $\bf \nu_n(q)$ and Proof of Proposition \ref{*q0=}}
\setcounter{equation}{0}

\begin{lemma} \label{*TraceD}
\noindent For each $q\!\in\!\bH_+$ the following identity holds:
\begin{equation}
\label{TraceNonLinD} \sum_{n\ge 0}(\sigma_n(q)\!-\!\sigma_n^0\!-\!2\nh q_{2n+1})=0\,.
\end{equation}
where the series converges absolutely. 
\end{lemma}

\begin{proof} By asymptotics (\ref{SAsympt}), the series $\sum_{n\ge 0}
(\sigma_n(q)\!-\!\sigma_n^0\!-\!2\nh q_{2n+1})$ converges
absolutely. Due to Lemma
\ref{*AAnalyticity}\begin{footnote}{Here and below $\partial\xi(q)\big/\partial
q\!=\!\zeta(q)$ means that for any $v\!\in\!L^2$ the equation $(d_q\xi)(v)\!=\!\langle
v,\overline{\zeta} \rangle_{L^2}$ holds true. }\end{footnote},
$
\frac{\partial\sigma_n(q)}{\partial q(x)}=\psi_{n,D}^2(x,q)\,,
$
where $\psi_{n,D}$ is the $n$-th normalized (in $L^2(\R_+)$) eigenfunction of the operator
$T_D$\,. Therefore,
$$
\sigma_n(q)-\sigma_n^0 = \int_0^1\frac{d}{ds}\,\sigma_{n}(sq)ds=
\int_0^1\langle\psi_{n,D}^2(x,sq),q(x)\rangle_+ds\,.
$$
Recall that $\psi_{n,D}^2(x,0)\!=\!2(\psi_{2n+1}^0)^2(x)$\,. Then,
$$
\sigma_n(q)\!-\!\sigma_n^0\!-\!2\nh q_{2n+1}= \int_0^1\!
\left\langle\psi_{n,D}^2(x,sq)\!-\!\psi_{n,D}^2(x,0),q(x)\right\rangle_{\!L^2(\R_+\,,\,dx)}ds\,.
$$
The standard perturbation theory (e.g., see \cite{Ka}) yields
$$
\frac{\partial\psi_{n,D}(x,q)}{\partial q(y)}= \sum_{m:m\ne n}
\frac{\psi_{n,D}(y,q)\psi_{m,D}(y,q)}{\sigma_{n}(q)\!-\!\sigma_{m}(q)}\,\psi_{m,D}(x,q)\,.
$$
Hence,
$$
\left\langle\psi_{n,D}^2(x,sq)\!-\!(\psi_{n,D}^0)^2(x),q(x)\right\rangle_{\!L^2(\R_+\,,\,dx)}=
\int_0^s\biggl\langle\frac{d}{dt}\,\psi_{n,D}^2(x,tq),q(x)\biggr\rangle_{\!L^2(\R_+\,,\,dx)}dt
$$
$$
= \int_0^s \biggl\langle 2\psi_{n,D}(x,tq)\cdot\biggl\langle\sum_{m:m\ne n}
\frac{\psi_{n,D}(y,tq)\psi_{m,D}(y,tq)}{\sigma_{n}(q)\!-\!\sigma_{m}(q)}\,\psi_{m,D}(x,tq),
q(y)\biggr\rangle_{\!L^2(dy)}\,, q(x) \biggr\rangle_{\!L^2(dx)}dt\,.
$$
This gives
$$
\sum_{n\ge 0} (\sigma_n(q)\!-\!\sigma_n^0\!-\!2\nh q_{2n+1})= 2\int_0^1 ds\int_0^s \sum_{n\ge
0} \sum_{m:m\ne n} \frac{\left\langle(\psi_{n,D}\psi_{m,D})(tq),q\right\rangle_+^2}
{\sigma_{n}(tq)\!-\!\sigma_{m}(tq)}\,dt\,.
$$
Let
$$
S_k=\sum_{n=0}^k\sum_{m:m\ne n}
\frac{\left\langle(\psi_{n,D}\psi_{m,D})(tq),q\right\rangle_+^2}
{\sigma_{n}(tq)\!-\!\sigma_{m}(tq)} = \sum_{n=0}^k\sum_{m=k+1}^{+\infty}
\frac{\left\langle(\psi_{n,D}\psi_{m,D})(tq),q\right\rangle_+^2}
{\sigma_{n}(tq)\!-\!\sigma_{m}(tq)}\,.
$$
Due to Lemma \ref{*xALemma}, for some absolute constant $\varepsilon\!>\!0$ we have
$$
\langle (\psi_{n,D}\psi_{m,D})(tq),q \rangle_+^2=\cases{O(n^{-\frac{1}{2}}\,m^{-\frac{1}{2}})&
for all $n,m\!\ge\!0$\,,\cr O(n^{-\frac{1}{2}-\frac{\varepsilon}{2}}m^{-\frac{1}{2}}), & if
$m\!\ge\!n\!+\!n^{\frac{1}{2}+\ve}$.}
$$
Using the simple estimate $|\sigma_n(tq)-\sigma_m(tq)|^{-1}\!=\!O(|n-m|^{-1})$ and technical
Lemma \ref{*ALemma}, we obtain $S_k\!\to\!0$ as $k\!\to\!\infty$\,, i.e. $\sum_{n\ge
0}(\sigma_n(q)\!-\!\sigma_n^0\!-\!2\nh q_{2n+1})=0$\,.
\end{proof}

\begin{proof}[{\bf Proof of Proposition \ref{*q0=}.}] Repeating
the proof of Lemma \ref{*TraceD}, we obtain  
$$
\sum_{n\ge 0} (\lambda_n(q)\!-\!\lambda_n^0\!-\!2\nh q_n) = 0,
$$
for $q\!\in\!\bH_{even}$\,.
Recall that $\lambda_{2n+1}(q)\!=\!\sigma_n(q)$ and $\lambda_{2n+1}^0\!=\!\sigma_n^0$\,. Using
(\ref{TraceNonLinD}), we get
$$
\sum_{n\ge 0} (-1)^n(\lambda_n(q)\!-\!\lambda_n^0\!-\!2\nh q_n) = \sum_{n\ge 0}
(\lambda_n(q)\!-\!\lambda_n^0\!-\!2\nh q_n) - 2\sum_{n\ge 0}
(\sigma_n(q)\!-\!\sigma_n^0\!-\!2\nh q_{2n+1})=0\,.
$$
Due to Lemma \ref{*NcNhasFG}, Propsition \ref{*H0Emb} (iv)  this yields
$$
\sum_{n\ge 0} (-1)^n(\lambda_n(q)\!-\!\lambda_n^0) = 2\sum_{n\ge 0} (-1)^n\nh q_n =
2\cdot\frac{(Fq)(z)}{\sqrt{1\!-\!z}}\Big|_{z=-1}=\frac{q(0)}{2}\,.
$$
Hence, $q(0)\!=\!2\sum_{n\ge 0}(-1)^n(\lambda_n(q)\!-\!\lambda_n^0)$,
which gives (\ref{*q0=}).
\end{proof}

Recall that each sequence $\{\sigma_n(q)\!-\!\sigma_n^0\}, q\!\in\!\bH_+$ belongs to the set
$\cS_D\!\subset\!\cH$ given by (\ref{SDDef}).
\begin{theorem}
\label{*RAsympt} (i) Each function
$r_n(q)=\nu_n(q)-\nu_n^0+\frac{q(0)}{2(2n\!+\!1)}, q\!\in\!\bH_+$ satisfies
\begin{equation}
\label{RAsympt} r_n(q)= -2\nt q_n +
R_n\left(\mu\right)+
\el2_{\frac{3}{4}+\delta}(n)\,,\ \ \mu=\{\mu_n\}_0^\infty,\ \
\mu_n=\sigma_m\!-\!\sigma_m^0
\end{equation}
uniformly on bounded subsets of $\bH_+$,where $\delta\!>\!0$ is some absolute constant  and  
\begin{equation}
\label{RDef}
\begin{array}{rl} \displaystyle
R_n \left(\mu\right)\, = & \displaystyle -2\log\,\biggl[
\frac{\psi_+^0(0,\sigma_n^0\!+\!\mu_n)}{\dot{\kappa}_{2n+1}\cdot \mu_n} \prod_{m:m\ne n}
\frac{(\sigma_n^0\!+\!\mu_n)-(\sigma_m^0\!+\!\mu_m)}{(\sigma_n^0\!+\!\mu_n)-\sigma_m^0}
\biggr] \cr & \displaystyle +\,\frac{2\dot{\kappa}'_{2n+1}}{\kappa'_{2n+1}}\cdot\mu_n -
\sum_{m\ge{0}}\frac{\mu_m}{2(n\!-\!m)\!+\!1}\,,\qquad n\!\ge\!0\,.
\end{array}
\end{equation}
(ii) For each $\{\mu_m\}_{m=0}^\infty\!\in\!\cS_D$ the sequence $\{R_n\}_{n=0}^\infty$ belongs
to the space $\el2_{3/4}$\,. Moreover, the mapping $\cR:\cS_D\to\el2_{3/4}$ given by $\{\mu_m\}_{m=0}^\infty \mapsto \{R_n\}_{n=0}^\infty$\,, is locally bounded. \remark {\rm Note
that $\nt q_n \!=\! \el2_{3/4}(n)$ due to Theorem \ref{*TildeQ} (ii). Therefore,
$\{r_{n}(q)\}_{n=0}^\infty \in\el2_{3/4}$\,.}
\end{theorem}

\begin{proof}
(i) Let $\sigma_n\!=\!\sigma_n(q)$ and $\mu_n\!=\!\sigma_n(q)\!-\!\sigma_n^0$\,,
$n\!\ge\!0$\,. Lemma \ref{*NuIdent} yields
$$
\frac{\nu_n(q)-\nu_n^0}{2} = \log \biggl[
\frac{\psi'_+(0,\sigma_n,q)}{\dot{\psi}_+(0,\sigma_n,q)}\cdot
\frac{\dot{\kappa}_{2n+1}}{\kappa'_{2n+1}} \biggr] \,.
$$
Using Theorem \ref{*BasicAsympt} (ii) and Theorem \ref{*TildeQ} (i), we obtain
$$
\log  \frac{\psi'_+(0,\sigma_n,q)}{\kappa'_{2n+1}} =
\frac{\dot{\kappa}'_{2n+1}}{\kappa'_{2n+1}}\,\mu_{n} - \frac{q(0)}{4(2n\!+\!1)}- \nt q_n -
\sum_{m\ge 0} \frac{\nh q_{2m+1}}{2(n\!-\!m)+1} + \el2_{\frac{3}{4}+\delta}(n)\,.
$$
Furthermore, identity (\ref{NuIdent2}) gives
$$
\log \frac{\dot{\psi}_+(0,\sigma_n,q)}{\dot{\kappa}_{2n+1}} = \log\,\biggl[
\frac{\psi_+^0(0,\sigma_n)}{\dot{\kappa}_{2n+1}\cdot \mu_n} \prod_{m:m\ne n}
\frac{\sigma_n-\sigma_m}{\sigma_n-\sigma_m^0} \biggr]
$$
Hence,
$$
r_n(q) = \nu_n(q)- \nu_n^0+ \frac{q(0)}{2(2n\!+\!1)} = -2\nt q_n +
R_n(\mu)+ \el2_{\frac{3}{4}+\delta}(n) + h_n\,,
$$
where
$$
h_n= \sum_{m\ge 0} \frac{\mu_m- 2\nh q_{2m+1}}{2(n\!-\!m)+1}\,,\quad n\!\ge\!0\,.
$$

In order to prove that $h_n\!=\!\el2_{\frac{3}{4}+\delta}(n)$, we note that identity
(\ref{SqrtDef}) yields
$$
h(z)\equiv\sum_{n\ge 0}h_nz^n\equiv
\frac{\pi}{2}P_+\biggl[\frac{g(\zeta)}{\sqrt{-\zeta}}\biggr]\,,\quad {\rm where}\quad
g(z)\equiv\sum_{m\ge 0}(\mu_{m}\!-\!2\nh q_{2m+1})z^m
$$
Due to asymptotics (\ref{SAsympt}) and identity (\ref{TraceNonLinD}), we have
$g\!\in\!H^2_{\frac{3}{4}+\delta}$ and $g(1)\!=\!0$\,. Hence\begin{footnote}{Here and below
$W^2_r(\T)$ is the Sobolev space on the unit circle
$\T\!=\!\{\zeta\!\in\!\C:|\zeta|\!=\!1\}$.}\end{footnote},
$$
\frac{g(\zeta)}{\sqrt{-\zeta}}\in W^2_{\frac{3}{4}+\delta}(\T)\,,\qquad {\rm and\ so}\qquad
P_+\biggl[\frac{g(\zeta)}{\sqrt{-\zeta}}\biggr]\in H^2_{\frac{3}{4}+\delta}\,.
$$
Thus, $\sum_{n\ge 0}h_nz^n\in H^2_{\frac{3}{4}+\delta}$\,, i.e.
$\{h_n\}_{n=0}^\infty\!\in\!\el2_{\frac{3}{4}+\delta}$\,.

\noindent (ii) Let $\{\mu_{m}\}_{m=0}^\infty\!\in\!\cS_D$\,. We rewrite (\ref{RDef}) in the
form $R_n=R_n^{(1)}+R_n^{(2)}+R_n^{(3)}$\,, where
$$
R_n^{(1)}= -2\log\,\biggl[\frac{\psi_+^0(0,\sigma_n^0\!+\!\mu_n)}
{\dot{\kappa}_{2n+1}\cdot\mu_n}\biggr]+
\frac{2\dot{\kappa}'_{2n+1}}{\kappa'_{2n+1}}\cdot\mu_n\,,
$$
$$
R_n^{(2)}= -\sum_{m:m\ne n} \left(2\log\, \biggl[1\!-\!\frac{\mu_m}{4(n\!-\!m)+\mu_n}\,\biggr]
+ \frac{\mu_m}{2(n\!-\!m)} \right)\,,
$$
$$
R_n^{(3)}= \frac{1}{2}\sum_{m:m\ne n} \frac{\mu_m}{n\!-\!m} - \sum_{m\ge 0}
\frac{\mu_m}{2(n\!-\!m)\!+\!1}\,.
$$
In the following Lemmas \ref{*xL51}\,--\,\ref{*xL53} we will analyze these terms separately.
Recall that Proposition \ref{*H0Emb} and Lemma \ref{*cH0Asymp} give $\mu_n\!=\!O(n^{-1/2})$ as
$n\!\to\!\infty$ and
\begin{equation}
\label{x54} \mu_n\!=\!v\cdot (n\!+\!1)^{-\frac{1}{2}}+ \el2_{1/4}(n)\,,
\end{equation}
where $v\!\in\!\R$ is some constant.

\begin{lemma}
\label{*xL51} The asymptotics $\displaystyle R_n^{(1)}= \frac{\pi^2v^2}{48}\,(n\!+\!1)^{-1} +
\el2_{3/4}(n)$ holds true.
\end{lemma}
\begin{proof}
Due to $\psi_+^0(0,\sigma_n^0)\!=\!0$, $\mu_n\!=\!O(n^{-1/2})$ and the estimates from
Corollary \ref{*ADotEstim}, we have
$$
\frac{\psi_+^0(0,\lambda_n^0\!+\!\mu_n)} {\dot{\kappa}_{2n+1}\cdot\mu_n}= 1
+\frac{\ddot\kappa_{2n+1}}{2\dot\kappa_{2n+1}}\,\mu_n +
\frac{\dddot\kappa_{2n+1}}{6\dot\kappa_{2n+1}}\,\mu_n^2+O(n^{-\frac{3}{2}}\log^4 n)\,.
$$
Therefore,
$$
R_{n}^{(1)}= \biggl(\frac{2\dot{\kappa}'_{2n+1}}{\kappa'_{2n+1}} -
\frac{\ddot\kappa_{2n+1}}{\dot\kappa_{2n+1}})\mu_{n} -2
\biggl(\frac{\dddot\kappa_{2n+1}}{6\dot\kappa_{2n+1}}-
\frac{(\ddot\kappa_{2n+1})^2}{8(\dot\kappa_{2n+1})^2}\biggr)\mu_n^2 + O(n^{-\frac{3}{2}}\log^4
n)
$$
$$
= \frac{\pi^2}{48}\,\mu_n^2 + O(n^{-\frac{3}{2}}\log^4 n)=
\frac{\pi^2v^2}{48}\,(n\!+\!1)^{-1}+\el2_{3/4}(n)\,,
$$
where we have used Lemma \ref{*LKappaSpecAsymp} and (\ref{x54}).
\end{proof}
\begin{lemma}
\label{*xL52} The asymptotics $\displaystyle R_n^{(2)}= -\frac{\pi^2v^2}{48}\,(n\!+\!1)^{-1} +
\el2_{3/4}(n)$ holds true.
\end{lemma}
\begin{proof} For $m\!\ne\!n$ we have
$$
2\log \, \biggl[1\!-\!\frac{\mu_m}{4(n\!-\!m)+\mu_n}\,\biggr] + \frac{\mu_m}{2(n\!-\!m)}
$$
$$
=\frac{\mu_m\mu_n}{2(n\!-\!m)(4(n\!-\!m)+\mu_n)} - \frac{\mu_m^2}{(4(n\!-\!m)+\mu_n)^2}+
O\biggl(\frac{m^{-3/2}}{(n\!-\!m)^3}\biggr)
$$
$$
=\frac{\mu_m\mu_n}{8(n\!-\!m)^2}- \frac{\mu_m^2}{16(n\!-\!m)^2} +
O\biggl(\frac{m^{-3/2}+m^{-1/2}n^{-1}}{(n\!-\!m)^3}\biggr)\,.
$$
Therefore,
$$
16R_n^{(2)}= -2\mu_n\cdot\!\!\sum_{m:m\ne n} \frac{\mu_m}{(n\!-\!m)^2}\ + \sum_{m:m\ne n}
\frac{\mu_m^2}{(n\!-\!m)^2}\, + O(n^{-\frac{3}{2}})\,.
$$
Recall that $\mu_n\!=\!v(n\!+\!1)^{-1/2}+\el2_{1/4}(n)$ and
$\mu_n^2\!=\!v^2(n\!+\!1)^{-1}+\el2_{3/4}(n)$\,. Using simple technical Lemma
\ref{*ALaboutL2}, we deduce that
$$
16R_n^{(2)}= -2\biggl(\frac{v}{(n\!+\!1)^{1/2}}+\el2_{\frac{1}{4}}(n)\biggr)
\biggl(\frac{\pi^2v}{3(n\!+\!1)^{1/2}}+\el2_{\frac{1}{4}}(n)\biggr) +
\biggl(\frac{\pi^2v^2}{3(n\!+\!1)}+\el2_{\frac{3}{4}}(n)\biggr)+O(n^{-\frac{3}{2}})\,.
$$
This gives $48R_n^{(2)}= -{\pi^2v^2}{(n\!+\!1)^{-1}}+\el2_{3/4}(n)$\,.
\end{proof}

\begin{lemma} \label{*xL53} The asymptotics $R_n^{(3)}=\el2_{3/4}(n)$ holds true.
\end{lemma}
\begin{proof}
Note that the following identities are fulfilled in $L^2(\T)$:
$$
2\sum_{l=-\infty}^{+\infty}\frac{\zeta^l}{2l\!+\!1} = \frac{\pi}{\sqrt{-\zeta}}\,, \qquad
\sum_{l:l\ne 0}\frac{\zeta^l}{l}=-\log(-\zeta)\,,
$$
where the branches of $\sqrt{-\zeta}$ and $\log(-\zeta)$, $\zeta\!\in\!\T\setminus\{1\}$ are
such that $\sqrt{1}\!=\!1$ and $\log 1\!=\!0$\,. Then,
$$
\sum_{n\ge 0}R_n^{(3)}z^n \equiv
-\frac{1}{2}\,P_+\biggl[\biggl(\frac{\pi}{\sqrt{-\zeta}}+\log(-\zeta)\biggr) \cdot \sum_{n\ge
0} \mu_n\zeta^n\biggr]
$$
Since $\{\mu_n\}_{n=0}^\infty\!\in\!\cH$\,, we have
$$
\sum_{n\ge 0}\mu_nz^n\equiv\frac{F(z)}{\sqrt{1\!-\!z}}\,,\quad {\rm where} \quad
F\!\in\!H^2_{3/4}\,.
$$
Introduce the function
\begin{equation}
\label{x5gamma}
\gamma(\zeta)=\frac{\pi\big/\sqrt{-\zeta}+\log(-\zeta)}{\sqrt{1\!-\!\zeta}}\,,\quad
\zeta\!\in\!\T\,.
\end{equation}
It is clear that $\gamma\!\in\!C^\infty(\T\setminus\{1\})$\,. Note that
$\gamma(\zeta)\!\to\!0$ as $\zeta\!\to\!1\!\pm\!i0$. This yields $\gamma\!\in\!W^2_{3/4}(\T)$,
$\gamma F\!\in\!W^2_{3/4}(\T)$ and $P_+[\gamma F]\!\in\!H^2_{3/4}$\,. The last statement is
equivalent to $\{R_n^{(3)}\}_{n=0}^\infty\!\in\!\el2_{3/4}$\,.
\end{proof}

Lemmas \ref{*xL51}\,--\,\ref{*xL53} give
$R_n\!=\!R_n^{(1)}\!+\!R_n^{(2)}\!+\!R_n^{(3)}\!=\!\el2_{3/4}(n)$\,. Note that all estimates
are uniform on bounded subsets of $\cS_D$\,. The proof of Theorem \ref{*RAsympt} is finished.
\end{proof}

\section{\label{SectInvD}Proof of Theorem \ref{*H+DChar}}
\setcounter{equation}{0}

Introduce the mapping
$$
\Phi: q\mapsto \biggl(\{\mu_n(q)\}_{n=0}^\infty\,;q(0)\,; \{r_n(q)\}_{n=0}^\infty\biggr)\,,
$$
where $\mu_n(q)\!=\!\sigma_n(q)\!-\!\sigma_n^0$\,. Due to Theorems \ref{*HevenChar},
\ref{*RAsympt}, we have
$$
\Phi: \bH_+\to \cS_D\times\R\times\el2_{\frac{3}{4}}\,.
$$
Theorem \ref{*H+DChar} claims that $\Phi$ is a real-analytic isomorphism. The proof given
below consists of five steps: $\Phi$ is injective (Sect. \ref{SectInvD}.1); $\Phi$ is
real-analytic (Sect. \nolinebreak \ref{SectInvD}.2); the Frech\'et derivative $d_q\Phi$ is a
Fredholm operator for each $q\!\in\!\bH_+$ (Sect. \ref{SectInvD}.3); $d_q\Phi$ is invertible
for each $q\!\in\!\bH_+$\,, i.e. $\Phi$ is a local real-analytic isomorphism (Sect.
\ref{SectInvD}.4); $\Phi$ is surjective (Sect. \ref{SectInvD}.5).

\begin{proof}[{\bf \ref{SectInvD}.1. Uniqueness Theorem.}]$\phantom{x}$

\vskip 6pt

\noindent Let $\left(\{\mu_n(p)\}_{n=0}^\infty\,;p(0)\,;
\{r_n(p)\}_{n=0}^\infty\right)\!=\!\left(\{\mu_n(q)\}_{n=0}^\infty\,;q(0)\,;
\{r_n(q)\}_{n=0}^\infty\right)$ for some $p,q\!\in\!\bH_+$. By definitions of $\mu_n$ and
$r_n$\,, it is equivalent to
$$
\sigma_n=\sigma_n(p)=\sigma_n(q)\quad {\rm and} \quad \nu_n=\nu_n(p)=\nu_n(q)\quad {\rm for \
all\ \ } n\!\ge\!0\,.
$$
Using Lemma \ref{*NuIdent}, we obtain
$$
\psi'_+(0,\sigma_n,p)=\psi'_+(0,\sigma_n,q)\quad {\rm for \ all \ \ } n\!\ge\!0\,.
$$
The rest of the proof is standard (see also \cite{CKK1}). Recall that $\varphi(x,\l,q)$ is the
solution of Eq. (\ref{PertEq}) such that $\vp(0,\l,q)\!=\!0$\,, $\vp'(0,\l,q)\!=\!1$\,.
Introduce the functions
$$
f_1(\l;x,q,p)= \frac{F_1(\l;x,q,p)}{\psi_+(0,\l,q)}\,,\ \ F_1(\l;x,q,p)=
\psi_+(x,\l,p)\varphi'(x,\l,q)-\varphi(x,\l,p)\psi'_+(x,\l,q)\,,
$$
$$
f_2(\l;x,q,p)= \frac{F_2(\l;x,q,p)}{\psi_+(0,\l,q)}\,,\ \ F_2(\l;x,q,p)=
\psi_+(x,\l,p)\varphi(x,\l,q) - \varphi(x,\l,p)\psi_+(x,\l,q)\,.
$$

Both $f_1$ and $f_2$ are entire with respect to $\l$ for each $x\!\in\!\R_+$\,. Indeed, all
roots $\sigma_n$\,, $n\!\ge\!0$\,, of the denominator $\psi_+(0,\cdot,q)$ are simple and all
these values are roots of the numerators $F_1$\,, $F_2$\,, since
$$
\frac{\psi_+(x,\sigma_n,p)}{\varphi(x,\sigma_n,p)} = \psi'_+(0,\sigma_n,p)=
\psi'_+(0,\sigma_n,q) = \frac{\psi_+(x,\sigma_n,q)}{\varphi(x,\sigma_n,q)}
$$
for all $x\!\in\!\R_+$ and $n\!\ge\!0$\,. Standard estimates (see Lemma \ref{*AEstim} and
asymptotics (\ref{EstimExact})) of $\varphi$ and $\psi_+$ give
$$
f_1(\l;x,p,q)=1\!+\!O(|\l|^{-\frac{1}{2}})\,,\quad f_2(\l;x,p,q)=O(|\l|^{-\frac{1}{2}})\,,
\quad |\l|\!=\!\l_{2k}^0\,,\ k\!\to\!\infty\,.
$$
Then, the maximum principle implies
$$
f_1(\l;x,p,q)= 1\,,\qquad f_2(\l;x,p,q)= 0\,, \quad\l\!\in\!\C\,.
$$
This yields $\vp(x,\l,p)=\vp(x,\l,q)$ and $\psi_+(x,\l,p)=\p_+(x,\l,q)$, i.e. $p=q$\,.
\end{proof}

\begin{proof}[{\bf \ref{SectInvD}.2. $\bf \Phi$ is a real-analytic mapping.}]
$\phantom{x}$

\vskip 6pt

\noindent Recall that for some $\delta\!>\!0$ the following asymptotics are fulfilled (see
\nolinebreak (\ref{SAsympt}) and (\ref{RAsympt})):
\begin{equation}
\label{MuSAsymptD} \mu_n(q)=2\nh q_{2n+1} + \el2_{\frac{3}{4}+\delta}(n)\,,\qquad r_n(q)=
-2\nt q_n + R_n\left(\{\mu_m(q)\}_{m=0}^\infty\right) + \el2_{\frac{3}{4}+\delta}(n)\,,
\end{equation}
where
$$
\cR:\cS_D\to \el2_{\frac{3}{4}}\,,\qquad \cR: \{\mu_m\}_{m=0}^\infty \mapsto
\{R_n\}_{n=0}^\infty\,,
$$
is a locally bounded mapping given by (\ref{RDef}). Let $\bH_{+\C}$ be the complexification of
$\bH_+$\,. Due to Lemma \ref{*AAnalyticity} (ii), for each $q\!\in\!\bH_+$ all functions
$\sigma_n(q)$ and $\psi'_+(0,\sigma_n(q),q)$\,, $n\!\ge\!0$\,, have analytic continuations
into some {complex} neighborhood of $q$\,. Moreover, due to Lemma \nolinebreak \ref{*NuIdent},
all functions $\nu_n(q)$ have analytic continuations into some complex neighborhood of $q$.
Therefore, for each {\it real} potential $q\!\in\!\bH_+$ all "coordinate functions"\
$\mu_n(q)$\,, $q(0)$\,, $r_n(q)$ of the mapping $\Phi$ have an analytic continuation into some
small {\it complex} neighborhood $Q$ of $q$.

Repeating the proof of (\ref{MuSAsymptD}) we obtain that these asymptotics hold true uniformly
on bounded subsets of ${Q}$\,. Let
$$
\Phi^{(0)}:q\mapsto\left(\{2\nh q_{2n+1}\}_{n=0}^\infty\,; q(0)\,; \{-2\nt q_n\}_{n=0}^\infty
\right)\,.
$$
Due to Theorem \ref{*TildeQ} (i), $\Phi^{(0)}$ is a linear isomorphism between $\bH_+$ and
$\cH\times\R\times\el2_{3/4}$\,. In particular, $\Phi^{(0)}$ is a real-analytic mapping.
Consider the difference
$$
\Phi-\Phi^{(0)}:q\mapsto \biggl(\{\mu_n(q)\!-\!2\nh q_{2n+1}\}_{n=0}^\infty\,; 0\,;
\{r_n(q)\!+\!2\nt q_n\}_{n=0}^\infty\biggr)\,,
$$
$$
\Phi-\Phi^{(0)}:\bH_+\to\el2_{\frac{3}{4}+\delta}\times\R\times\el2_{\frac{3}{4}}\,.
$$
All "coordinate functions"\ $\mu_n(q)\!-\!2\nh q_{2n+1}$\,, $r_n(q)\!+\!2\nt q_n$ are analytic
and $\Phi\!-\!\Phi^{(0)}$ is correctly defined and bounded in some small complex neighborhood
of each real potential (since (\ref{MuSAsymptD}) holds true uniformly on bounded subsets).
Then, $\Phi\!-\!\Phi^{(0)}$ is a real-analytic mapping from $\bH_+$ into
$\el2_{\frac{3}{4}+\delta}\!\times\!\R\!\times\!\el2_{\frac{3}{4}}$ and $\Phi$ is
real-analytic too, since
$\el2_{\frac{3}{4}+\delta}\!\subset\!\el2_{\frac{3}{4}}\!\subset\!\cH$\,.
\end{proof}

\begin{proof}[{\bf \ref{SectInvD}.3. The Frech\'et derivative
$\bf d_q\Phi$ is a Fredholm operator for each $q\!\in\!\bH_+$\,.}] $\phantom{x}$

\vskip 6pt

\noindent In other words, we will prove that $d_q\Phi$ is the sum of invertible and compact
operators. Let
$$
\Phi^{(1)}:q\mapsto \left(0\,;0\,;\cR(\{\mu_{m}(q)\}_{m=0}^\infty)\right)\quad {\rm and}\quad
\Phi^{(2)}=\Phi\!-\!\Phi^{(0)}\!-\!\Phi^{(1)}\,.
$$
Using the same arguments as above, we obtain that $\cR_D:\cS_D\to\el2_{3/4}$ is a
real-analytic mapping (since it is locally bounded in some small complex neighborhood of each
real point $\mu\!\in\!\cS_D$ and all "coordinate function"\ $R_n$ are analytic). Then,
$\Phi^{(1)}$ is real-analytic as a composition of real-analytic mappings. Theorems
\ref{*BasicAsympt}, \ref{*RAsympt} yield
$$
\Phi^{(2)}:\bH_+\to \el2_{\frac{3}{4}+\delta}\times\R\times\el2_{\frac{3}{4}+\delta}\,.
$$
Repeating above arguments again, we obtain that $\Phi^{(2)}$ is a real-analytic mapping too.

Fix some $q\!\in\!\bH_+$\,. The Frech\'et derivatives $d_q\Phi$\,, $d_q\Phi^{(j)}$ of the
analytic mappings $\Phi$, $\Phi^{(j)}$ at the point $q$ are bounded linear operators and
$$
d_q\Phi=d_q\Phi^{(0)}+d_q\Phi^{(1)}+d_q\Phi^{(2)}=
(\Phi^{(0)}\!+\!d_q\Phi^{(1)})+d_q\Phi^{(2)}\,.
$$
Note that the operator
$$
d_q\Phi^{(2)}:\bH_+\to \cH\times\R\times\el2_{\frac{3}{4}}
$$
is compact since it maps $\bH_+$ into
$\el2_{\frac{3}{4}+\delta}\times\R\times\el2_{\frac{3}{4}+\delta}$ and the embedding
$\el2_{\frac{3}{4}+\delta}\subset \el2_{\frac{3}{4}}$ is compact. In order to prove that
$\Phi^{(0)}\!+\!d_q\Phi^{(1)}$ is invertible, we introduce two linear operators
$$
A:\bH_+\to\cH\,,\ \ \ p\mapsto Ap=\{2\nh p_{2n+1}\}_{n=0}^\infty\,,\qquad
B:\bH_+\to\el2_{\frac{3}{4}}\,,\ \ \ p\mapsto Bp=\{-2\nt p_n\}_{n=0}^\infty\,.
$$
Recall that $\Phi^{(0)}p\!=\!(Ap\,;p(0)\,;Bp)$\,. The chain rule implies
$$
(d_q\Phi^{(1)})p=(0\,;0\,; (d_{\mu(q)}\cR)Ap)\,,
$$
where $d_{\mu(q)}\cR$ is the Frech\'et derivative of the mapping $\cR$ at the point
$\mu(q)\!=\!\{\mu_m(q)\}_{m=0}^\infty\!\in\!\cS_D$\,. Hence,
$\Phi^{(0)}\!+\!d_q\Phi^{(1)}\!=\!C\Phi^{(0)}$\,, where both operators $C$ and $C^{-1}$ given
by
$$
C^{\pm 1}:\cH\times\R\times\el2_{3/4}\to\cH\times\R\times\el2_{3/4}\,,\quad C:(h;t;r)\mapsto
\left(h;t;r\!\pm\!(d_{\mu(q)}\cR_D)h\right)
$$
are bounded. Recall that $(\Phi_D^{(0)})^{-1}$ is bounded due to Theorem \ref{*TildeQ} (ii).
Therefore, the operator $(\Phi^{(0)}\!+\!d_q\Phi^{(1)})^{-1}$ is bounded too.
\end{proof}

\begin{proof}[{\bf \ref{SectInvD}.4. $\bf\Phi$ is a local real-analytic isomorphism.}]
$\phantom{x}$

\vskip 6pt

\noindent By Fredholm's Theory, in order to prove that $(d_q\Phi)^{-1}$ is bounded, it is
sufficient to check that the range $\Ran d_q\Phi$ is dense:
\begin{equation}
\label{ClImPhi} \cH\times\R\times\el2_{3/4}=\overline{\Ran d_q\Phi}\,.
\end{equation}
Note that Lemma \ref{*AAnalyticity} (ii) gives
\begin{equation}
\label{GradIdent} \frac{\partial\sigma_n(q)}{\partial q(t)}= \psi_{n,D}^2(t,q)\,,\qquad
\frac{\partial \log[(-1)^n\psi'_+(0,\sigma_n(q),q)]}{\partial q(t)}=
-(\psi_{n,D}\chi_{n,D})(t,q)\,,
\end{equation}
where $\psi_{n,D}(\cdot,q)$ is the $n$-th normalized eigenfunction of $T_D$ and
$\chi_{n,D}(\cdot,q)$ is some special solution of Eq. (\ref{PertEq}) for $\l\!=\!\sigma_n(q)$
such that $\{\chi_{n,D}\,,\psi_{n,D}\}\!=\!1$\,. In particular,
$$
(\psi_{n,D}\chi_{n,D})(t,q)\sim t\,,\ \ t\!\to\!0\,,\qquad (\psi_{n,D}\chi_{n,D})(t,q)\sim
-t^{-1},\ \ t\!\to\!+\infty\,.
$$

Due to Lemma \ref{*AGradProd} (i), for each $q\in\bH_+$ the following standard identities are
fulfilled:
\begin{equation} \label{GradProd}
\begin{array}{c}\displaystyle
\left((\psi_{n,D}^2)'(q),\psi_{m,D}^2(q)\right)_{+}=0\,,\quad
\left((\psi_{n,D}^2)'(q),(\psi_{m,D}\chi_{2m+1})(q)\right)_{+}=
{\textstyle\frac{1}{2}}\,\delta_{mn}\,, \vphantom{\big|_{\big|}}\cr \displaystyle
\vphantom{\big|^{\big|}}
\left((\psi_{n,D}\chi_{n,D})'(q),(\psi_{m,D}\chi_{m,D})(q)\right)_{+}=0\,,\quad n,m\!\ge\!0.
\end{array}
\end{equation}
Note that $(\psi_{m,D}^2)'(\cdot,q)\!\in\!\bH_+$\,. Using (\ref{GradIdent}), (\ref{GradProd}),
we obtain
$$
\biggl(\frac{\partial\mu_n(q)}{\partial q}\,,(\psi_{m,D}^2)'(q)\biggr)_{\!+}\!\!=0\quad {\rm
and} \quad \biggl(\frac{\partial \log[(-1)^n\psi'_+(0,\sigma_n(q),q)]}{\partial
q}\,,(\psi_{m,D}^2)'(q)\biggr)_{\!+}\!\!= \frac{\delta_{nm}}{2}
$$
for all $n,m\!\ge\!0$\,. Due to Lemma \ref{*NuIdent}, this implies
$$
\biggl(\frac{\partial\dot{\psi}_+(0,\sigma_n(q),q)}{\partial
q}\,,(\psi_{m,D}^2)'(q)\biggr)_{\!+}\!\!=0 \qquad {\rm and} \qquad \biggl(\frac{\partial
\nu_n(q)}{\partial q}\,,(\psi_{m,D}^2)'(q)\biggr)_{\!+}\!\!= \delta_{nm}\,.
$$
The identity
$$
\biggl(\frac{\partial q(0)}{\partial q}\,,(\psi_{m,D}^2)'(q)\biggr)_{\!+}\!\!=
(\psi_{m,D}^2)'(0)=0
$$
gives
$$
\biggl(\frac{\partial r_n(q)}{\partial q}\,,(\psi_{m,D}^2)'(q)\biggr)_{\!+}\!\!=
\biggl(\frac{\partial \nu_n(q)}{\partial q}\,,(\psi_{m,D}^2)'(q)\biggr)_{\!+}\!\!=
\delta_{nm}\,.
$$
Thus,
$$
(d_q\Phi)\left((\psi_{m,D}^2)'(q)\right) = \left({\bf 0}\,; 0\,;{\bf e}_m\right)\,,
$$
where ${\bf 0}\!=\!(0,0,0,\dots)$\,, ${\bf e}_0\!=\! (1,0,0,\dots)$\,, ${\bf
e}_1\!=\!(0,1,0,\dots)$ and so on. Therefore,
\begin{equation}
\label{xClImPhi} \left\{({\bf 0};0)\right\}\times\el2_{3/4} =\{({\bf 0};0;{\bf c}):{\bf
c}\!\in\!\el2_{3/4}\} \subset\overline{\Ran d_q\Phi}\,.
\end{equation}

We come to the second component of $(d_q\Phi)\xi$\,, i.e. to the value $\xi(0)$\,. We consider
the lowest eigenvalue $\l_0(q)$ of the operator $T_N$ (with the same potential $q$ and the
Neumann boundary condition $\psi'(0)\!=\!0$) and the function
$$
\xi(t)= (\varphi\vartheta)'(t,\l_0(q),q)\,,\quad t\in \R\,,
$$
where $\vartheta(t)$ is the solution of $-\psi''\!+\!x^2\psi\!+\!q(x)\psi\!=\!\lambda\psi$
such that $\vartheta(0)\!=\!1$ and $\vartheta'(0)\!=\!1$\,. Note that $\xi(0)\!=\!1$\,.
Asymptotics (\ref{Psi+Asympt1}), (\ref{Chi+Asympt1}) give $\xi\!\in\!\bH_+$ since
$\vartheta(\cdot,\l_0(q),q)$ is proportional to $\psi_+(\cdot,\l_0(q),q)$\,. Moreover, using
$\varphi(0,\lambda_0(q),q)\!=\!\psi_{n,D}(0,q)\!=\!0$\,, we obtain
$$
\biggl(\frac{\partial\mu_n(q)}{\partial q}\,,\xi\biggr)_{\!+}\!\!=
\left(\psi_{n,D}^2(q)\,,\xi\right)_{+}=
\frac{-\{\psi_{n,D}\,,\varphi\}\{\psi_{n,D}\,,\vartheta\}(0,\l_0(q),q)}
{2(\sigma_n(q)\!-\!\l_0(q))}=0\,.
$$
Hence,
$$
(d_q\Phi)\xi = \left({\bf 0}\,; 1\,; (d_q r)\xi \right)\,,\quad  {\rm where}\quad (d_q r)\xi=
\biggl\{\biggl(\frac{\partial r_n(q)}{\partial q}\,,\xi\biggr)_{\!+}\biggr
\}_{n=0}^{\infty}\in\el2_{3/4}\,.
$$
Together with (\ref{xClImPhi}) this implies
$$
\{{\bf 0}\}\times\R\times\el2_{3/4}\subset\overline{\Ran d_q\Phi}\,.
$$

Furthermore, we consider the functions $-2(\psi_{m,D}\chi_{m,D})'(q)\in\bH_+$ (see asymptotics
(\ref{Psi+Asympt1}), (\ref{Chi+Asympt1})). Identities (\ref{GradIdent}), (\ref{GradProd}) and
$(\psi_{m,D}\chi_{m,D})'(0,q)\!=\!1$ give
$$
(d_q\Phi)\left(-2(\psi_{m,D}\chi_{m,D})'(q)\right) = \left({\bf e}_m\,;
-2\,;(d_qr)\left(-2(\psi_{m,D}\chi_{m,D})'(q)\right) \right)\,.
$$
Due to Proposition \ref{*H0Emb} (ii), the set of finite sequences is dense in $\cH_0$\,.
Therefore,
\begin{equation}
\label{xxClImPhi} \cH_0\times\R\times\el2_{3/4}\subset\overline{\Ran d_q\Phi}\,.
\end{equation}
In conclusion, we consider an arbitrary function $\zeta\!\in\!\bH_+$ such that
$\int_{\R_+}\!\zeta(t)dt\!\ne\!0$. Proposition \ref{*H0Emb} (i) implies
$$
(d_q\Phi)\zeta\notin\cH_0\times\R\times\el2_{3/4}\,.
$$
Together with (\ref{xxClImPhi}) this yields (\ref{ClImPhi}), since the codimension of $\cH_0$
in $\cH$ is equal to $1$.
\end{proof}

\begin{proof}[{\bf \ref{SectInvD}.4. $\Phi$ is surjective.}]
$\phantom{x}$

\vskip 6pt

\begin{lemma}
\label{*S2k+1Change} Let $q\!\in\!\bH_+$\,, $n\!\ge\!0$ and $t\!\in\!\R$\,. Denote
$$
q_n^t(x)=q(x)-2\frac{d^2}{dx^2}\log\eta_n^t(x,q)\,, \quad \eta_n^t(x,q)=
1+(e^t\!-\!1)\int_x^{+\infty}\!\!\psi_{n,D}^2(s,q)ds\,.
$$
Then $q_n^t\!\in\!\bH_+$ and
$$
\sigma_m(q_n^t)=\sigma_m(q)\,, \qquad \nu_m(q_n^t)=\nu_m(q)+t\delta_{mn}
$$
for all $m\!\ge\!0$\,. Moreover, $q_n^t(0)\!=\!q(0)$\,. \remark {\rm Therefore,
$r_m(q_n^t)\!=\!r_m(q)\!+\!t\delta_{nm}$ for all $m\!\ge\!0$\,.}
\end{lemma}
\begin{proof}
This Lemma is similar to \cite{CKK} Theorem 3.5 and can be proved by direct calculations using
the so-called Darboux transform of second-order differential equation (see also \cite{MT},
\cite{PT}) and Lemma \ref{*NuIdent}. Note that $\eta_n^t(x,q)\!=\!
e^t\!-\!(e^t\!-\!1)\int_0^x\!\!\psi_{n,D}^2(s,q)ds\!=\!e^t\!+\!O(x^3)$\,, $x\downarrow 0$\,.
This implies $q_n^t(0)=q(0)$\,.
\end{proof}

We consider an arbitrary spectral data $(h^*;u^*;c^*)\in\cS_D\times\R\times\el2_{3/4}$\,. Due
to Theorem \nolinebreak \ref{*HevenChar} and Proposition \ref{*q0=}, there exists a potential
$q^*\!\in\!\bH_+$ such that
$$
\mu_n(q^*)=h^*_n\ \ {\rm for\ all}\ \ n\!\ge\!0\quad {\rm and}\quad q^*(0)=u^*\,.
$$
This yields
$$
\left(h^*;u^*;r(q^*)\right)\in \Phi(\bH_+)\,,\quad {\rm where}\quad
r(q^*)\!=\!(r_0(q^*)\,,r_1(q^*)\,,\dots)\,.
$$
Due to Proposition \ref{*H0Emb} (ii), for each $\varepsilon\!>\!0$ there exist a finite
sequence $t_\varepsilon\!=\!(t_0\,,...\,,t_k\,,0\,,...)$ such that
$$
\left\|(c^*\!-t_\varepsilon)-r(q^*)\right\| = \left\|(c^*\!-r(q^*))-t_\varepsilon\right\| <
\varepsilon\,.
$$
Since $\Phi$ is a local isomorphism, for some $\varepsilon\!>\!0$ we have
$$
\left(h^*;u^*;c^*\!-\!t_\varepsilon\right)=
\left(h^*;u^*;(c^*_0\!-\!t_0\,,\dots,c^*_k\!-\!t_k\,, c^*_{k+1}\,,c^*_{k+2}\,,\dots)\right)\in
\Phi(\bH_+)\,.
$$
It means that $\left(h^*;u^*;c^*\!-\!t_\varepsilon\right)\!=\!\Phi(q_{k+1})$ for some
$q_{k+1}\!\in\!\bH_+$. Using Lemma \ref{*S2k+1Change} step by step, we construct the sequence
of potentials
$$
q_{j}=(q_{j+1})_{j}^{t_j}\!\in\!\bH_+\,,\quad j\!=\!k,k\!-\!1,\dots,1,0\,,
$$
such that
$$
\Phi(q_j)=\left(h^*;u^*; (c^*_0\!-\!t_0\,,\dots,c^*_{j-1}\!-\!t_{j-1}\,, c^*_j\,,
c^*_{j+1}\,,\dots)\right)\,.
$$
Then, $\Phi(q_0)\!=\!(h^*;u^*;c^*)$\,.
\end{proof}




\renewcommand{\thesection}{A}
\section{Appendix}
\setcounter{equation}{0}
\renewcommand{\theequation}{A.\arabic{equation}}

Here we collect some technical results from \cite{CKK1}, \cite{CKK} which are essentially used
above.

\vskip 6pt

\noindent {\bf A.1 The unperturbed equation. }

\vskip 6pt

\noindent For each $\l\!\in\!\C$ the equation $-\psi''\!+\!x^2\psi\!=\!\lambda\psi$ has the
solution $\p_{+}^0(x,\l)=D_{\frac{\l-1}{2}}(\sqrt{2}x)$\,, where $D_{\mu}(x)$ is the Weber
function (or the parabolic cylinder function, see \cite{B}). For each $x$ the functions
$\p_{+}^0(x,\cdot)$ and $(\p_{+}^0)'(x,\cdot)$ are entire and the following asymptotics are
fulfilled:
\begin{equation}
\label{Psi+0Asympt}
\begin{array}{c}\displaystyle
\p_{+}^0(x,\l)=(\sqrt{2}x)^{\frac{\l-1}{2}}e^{-\frac{x^2}{2}} \left(1+O(x^{-2})\right)\,,\quad
x\to +\infty\,, \vphantom{\frac{1}{\sqrt{2}}}\cr \displaystyle \label{Psi+0DerAsympt}
(\p_{+}^0)'(x,\l)=-\frac{1}{\sqrt{2}}(\sqrt{2}x)^{\frac{\l+1}{2}}e^{-\frac{x^2}{2}}
\left(1+O(x^{-2})\right)\,,\quad x\to +\infty\,,
\end{array}
\end{equation}
uniformly with respect to $\l$ on bounded domains. Note that (see \cite{B})
\begin{equation}
\label{Psi+00lambda}
\begin{array}{c}\displaystyle\psi_+^0(0,\l)=
D_{\frac{\l-1}{2}}(0)=2^{\frac{\l-1}{4}}\,\frac{\Gamma(\frac{1}{2})}{\Gamma(\frac{3-\l}{4})}=
\cos\frac{(\l\!-\!1)\pi}{4}\,\cdot\,
\frac{2^\frac{\l-1}{4}}{\sqrt{\pi}}\,\,\Gamma\!\biggl(\frac{\l\!+\!1}{4}\biggr)\,,\cr
\displaystyle (\psi_+^0)'(0,\l)=\sqrt{2}D'_{\frac{\l-1}{2}}(0)=
2^{\frac{\l-1}{4}}\,\frac{\Gamma(-\frac{1}{2})}{\Gamma(\frac{1-\l}{4})}=
\sin\frac{(\l\!-\!1)\pi}{4}\,\cdot\,
\frac{2^\frac{\l+3}{4}}{\sqrt{\pi}}\,\,\Gamma\!\biggl(\frac{\l\!+\!3}{4}\biggr)\,.
\end{array}
\end{equation}
Let $J^0(x,t;\l)$ be the solution of $-\psi''\!+\!x^2\psi\!=\!\lambda\psi$ such that
$J^0(t,t;\l)= 0$, $(J^0)'_x(t,t;\l)= 1$\,. Then
\begin{equation} \label{TasJ}
J^0(0,t;\l)= -\varphi^0(t,\l)=-\varphi(t,\l,0)\,,\qquad (J^0)'_x(0,t;\l)= \vt^0(t,\l)=\vt(t,\l,0)\,.
\end{equation}

In order to estimate $\psi_+^0$ and $J^0$, we introduce real-valued functions
\begin{equation}
\label{ASRdef}
\begin{array}{c}
\displaystyle a(\l)=
\biggl|\frac{\l}{2e}\biggr|^{\frac{\Re\l}{4}}\!e^{\frac{\pi-\phi}{4}\Im\l}, \quad
\l=|\l|e^{i\phi},\ \ \phi\in[0,2\pi)\,, \vphantom{\int_{\Big|}}\cr \displaystyle \r(x,\l)=
1+|\l|^{1/12}+|x^2-\l|^{1/4}\,,\qquad \s(x,\l)= \Re\int_0^x\!\!\sqrt{y^2-\l}\,dy\,,\ \
x\!\ge\!0\,,
\end{array}
\end{equation}
where $\sqrt{y^2\!-\!\l}=y\!+\!o(1)$ as $y\!\to\!+\infty$ (it is equivalent to
$\Re\sqrt{y^2\!-\!\l}\!\ge\!0$\,, if $y\!\ge\!0$).

\begin{lemma} \label{*AUnpertEstim}
For all $(x,t,\l)\in\R_+\times\R_+\times\C$ the following estimates are fulfilled:
\begin{equation}\label{Psi0Estim}
|\p_{+}^0(x,\l)|\le C_0 a(\l)\cdot \frac{e^{-\s(x,\l)}}{\r(x,\l)}\,, \qquad\qquad\
|(\p_{+}^0)'(x,\l)|\le C_0 a(\l)\cdot \r(x,\l)e^{-\s(x,\l)}\,,
\end{equation}
$$
|J^0(x,t;\l)|\le \frac{C_1}{\r(x,\l)\r(t,\l)}\, e^{|\s(x,\l)-\s(t,\l)|}\,, \quad
|(J^0)'_x(x,t;\l)|\le C_1\frac{\r(x,\l)}{\r(t,\l)}\, e^{|\s(x,\l)-\s(t,\l)|}\,,
$$
where $C_0,\ C_1$ are some absolute constants.
\end{lemma}
\begin{proof} See Lemmas 2.1 and 2.3 \cite{CKK1}. Note that the proof is based on
the result of \cite{O}. \end{proof}

\remark If $x\!=\!0$ and $|\lambda|\!\ge\!1$, then\begin{footnote}{Here and below $f\asymp g$
means that $C_1|f|\le|g|\le C_2|f|$ for some absolute constants $C_1,C_2>0$.}\end{footnote}
$\s(0,\lambda)\!=\!0$ and $\rho(0,\lambda)\!\asymp |\lambda|^{1/4}$. It follows from
identities (\ref{Psi+00lambda}) and routine calculations that
\begin{equation}\label{EstimExact}
\begin{array}{ll}\displaystyle
|\psi_+^0(0,\lambda)|\asymp |\lambda|^{-{1}/{4}}a(\lambda)\,, & {\rm if}\ \ |\lambda|=k\ne
\lambda_{2n+1}^0\,,\vphantom{\big|_{\big|}}\cr\displaystyle |(\psi_+^0)'(0,\lambda)|\asymp
|\lambda|^{{1}/{4}}a(\lambda)\,,& {\rm if}\ \ |\lambda|=k\ne \lambda_{2n}^0\,,
\vphantom{\big|^{\big|}}\end{array}\ \ k,n\!\in\!\N\,.
\end{equation}
In other words, the estimates (\ref{Psi0Estim}) of $|\psi_+(0,\lambda)|$ and
$|\psi'_+(0,\lambda)|$ are exact on these contours.

\vskip 6pt

\noindent {\bf A.2 The perturbed equation. }

\vskip 6pt

\noindent The solutions $\psi_+\,,\vartheta\,,\varphi$ of the perturbed equation
$-\psi''\!+\!x^2\psi\!+\!q(x)\psi\!=\!\lambda\psi$\,, $\lambda\!\in\!\C$\,, can be constructed
by iterations:
\begin{equation}
\label{Psi+asSum} \p_{+}(x,\l,q)={\mathop{\sum}\limits_{n\ge 0}} \p_{+}^{(n)}(x,\l,q)\,, \quad
\p_{+}^{(n+1)}(x,\l,q)= -\int_x^{+\infty}\!\!J^0(x,t;\l)\p_{+}^{(n)}(t,\l,q)q(t)dt\,,
\end{equation}
\begin{equation}
\label{T12asSum} \vt_{1,2}(x,\l,q)={\mathop{\sum}\limits_{n\ge
0}}\vt_{1,2}^{(n)}(x,\l,q),\qquad \vt_{1,2}^{(n+1)}(x,\l,q)=
\int_0^x\!J^0(x,t;\l)\vt_{1,2}^{(n)}(t,\l,q)q(t)dt\,,
\end{equation}
where we use the notations $\vt_1\!=\!\vt$ and $\vt_2\!=\!\varphi$ for short and $\vt_1^{(0)}=\vt^{(0)}, \vp_2^{(0)}=\vp^{(0)}$, see \er{TasJ}.

Introduce functions
$$
\beta_{+}(x,\l,q)= C_1\int_x^{+\infty}\!\frac{|q(t)|dt}{\r^2(t,\l)}\,, \qquad
\beta_{0}(x,\l,q)= C_1\int_0^x\frac{|q(t)|dt}{\r^2(t,\l)}dt\,.
$$
It is easy to see (\cite{CKK} Lemma \nolinebreak 5.5) that
\begin{equation}
\label{BetaEstim} \beta(\lambda,q)=\beta_{+}(x,\l,q)\!+\!\beta_{0}(x,\l,q)=
C_1\int_0^{+\infty}\frac{|q(t)|dt}{\r^2(t,\l)}=O(|\lambda|^{-1/2}\|q\|_{\bH_+})\,.
\end{equation}

\begin{lemma} \label{*AEstim}
For all $(x,\l,q)\in\R_+\!\times\C\times\bH_{+\C}$ the following estimates are fulfilled:
$$
|\p_{+}^{(n)}(x,\l,q)|\le C_0 a(\l)\frac{e^{-\s(x,\l)}}{\r(x,\l)}\cdot
\frac{\beta_{+}^n(x,\l,q)}{n!}\,,
$$
$$
|\vt_j^{(n)}(x,\l,q)|\le \frac{2 C_1}{(1\!+\!|\l|^{1/4})^{2j-3}}
\cdot\frac{e^{\s(x,\l)}}{\r(x,\l)}\cdot \frac{\beta_0^n(x,\l,q)}{n!}\,,\quad j\!=\!1,2\,.
$$
In particular, series (\ref{Psi+asSum}), (\ref{T12asSum}) converge uniformly on bounded
subsets of $\R_+\times\C\times\bH_{+\C}$\,. Moreover, the similar estimates with $\rho(x,\l)$
instead of $\frac{1}{\rho(x,\l)}$ in right-hand sides hold true for the values
$|(\p_{\pm}^{(n)})'(x,\l,q)|$ and $|(\vt_j^{(n)})'(x,\l,q)|$\,.
\end{lemma}
\begin{proof} See \cite{CKK1} Lemma 3.1 and \cite{CKK} Lemmas 5.2, 5.3.\end{proof}

\begin{corollary}\label{*ADotEstim}
For all $(\l,q)\in\C\times\bH_{+\C}$\,, $n,m\!\ge\!0$ and some absolute constant $C\!>\!0$ the
following estimates are fulfilled:
\begin{equation}\label{APsiNMEstim}
\begin{array}{c}\displaystyle
\biggl|\frac{\partial^m\p_{+}^{(n)}(0,\l,q)}{\partial\lambda^m}\biggr| \le
\frac{m!C^{n+m+1}\|q\|_{\bH_+}^n}{n!}\cdot \frac{\log^m(|\lambda|\!+\!2)\cdot
a(\l)}{(|\lambda|+1)^{\frac{n}{2}+\frac{1}{4}}}\,, \vphantom{\Big|_{\Big|}}\cr\displaystyle
\biggl|\frac{\partial^m(\p_{+}^{(n)})'(0,\l,q)}{\partial\lambda^m}\biggr| \le
\frac{m!C^{n+m+1}\|q\|_{\bH_+}^n}{n!}\cdot \frac{\log^m(|\lambda|\!+\!2)\cdot
a(\l)}{(|\lambda|+1)^{\frac{n}{2}-\frac{1}{4}}}\,.\vphantom{\Big|^{\Big|}}
\end{array}
\end{equation}
\end{corollary}
\begin{proof} Note that $\sigma(0,\lambda)\!=\!0$ and
$\rho(0,\lambda)\!\asymp\!1\!+\!|\lambda|^{1/4}$\,. Hence, Lemma \ref{*AEstim} and
(\ref{BetaEstim}) give (\ref{APsiNMEstim}) for $m\!=\!0$\,. Recall that
$\psi_+^{(n)}(0,\lambda,q)$, $(\psi_+^{(n)})'(0,\lambda,q)$ are entire functions. Therefore,
the simple estimate
$$
a(\lambda(\phi))\!=\!O(a(\lambda))\,,\quad {\rm if}\quad
\lambda(\phi)\!=\!\lambda+e^{i\phi}\log^{-1}(|\lambda|\!+\!2)\,,
$$
and the integration over the contour $\lambda(\phi)$\,, $\phi\!\in\![0,2\pi]$\,, imply
(\ref{APsiNMEstim}) in the case $m\!>\!0$\,.
\end{proof}

Let
$$
\widetilde{\beta}_{+}(x,q)=\frac{1}{x^2\!+\!1}\,+
\int_x^{+\infty}\biggl|\frac{q(t)}{t}\biggr|\,dt\,.
$$
The following asymptotics as $x\!\to\!+\infty$ are fulfilled uniformly on bounded subsets of
$\C\times\bH_{+\C}$ (see \cite{CKK} p.139 and p.169):
\begin{equation}
\label{Psi+Asympt1}
\begin{array}{c}\displaystyle
\p_{+}(x,\l,q)=(\sqrt{2}x)^{\frac{\l-1}{2}}e^{-\frac{x^2}{2}}
(1\!+\!O(\widetilde{\beta}_+(x,q)))\,, \vphantom{\frac{1}{\sqrt{2}}}\cr \displaystyle
\p'_+(x,\l,q)=-\frac{1}{\sqrt{2}}(\sqrt{2}x)^{\frac{\l+1}{2}}e^{-\frac{x^2}{2}}
(1\!+\!O(\widetilde{\beta}_+(x,q)))\,.
\end{array}
\end{equation}

\noindent Moreover, if $\chi_{+}(x,\l,q)$ is a solution of
$-\psi''\!+\!x^2\psi\!+\!q(x)\psi\!=\!\lambda\psi$ such that $k=\{\chi_{+},\p_{+}\}\ne 0$\,,
then
\begin{equation}
\label{Chi+Asympt1}
\begin{array}{c}\displaystyle
\chi_{+}(x,\l,q)=-\frac{k}{\sqrt{2}} (\sqrt{2}x)^{\frac{-\l-1}{2}}e^{\frac{x^2}{2}}
(1\!+\!O(\widetilde{\beta}_{+}(x,q)))\,, \cr
\displaystyle\vphantom{\frac{1^|}{\sqrt{2}}}\chi'_{+}(x,\l,q)\!=\!-\frac{k}{2}
(\sqrt{2}x)^{\frac{-\l+1}{2}}e^{\frac{x^2}{2}} (1\!+\!O(\widetilde{\beta}_{+}(x,q)))\,.
\end{array}
\end{equation}
\remark If $q\!\in\!\bH_{+\C}$\,, then (\ref{Psi+Asympt1}), (\ref{Chi+Asympt1}) give
$(\psi_+\chi_+)'\!\in\!\bH_{+\C}$ (see \cite{CKK} p. 172).


\vskip 6pt

\noindent{\bf A.3 Analyticity of spectral data and its gradients.}

\vskip 6pt

\noindent Recall that $\bH_{+\C}$ is the complexification of the space $\bH_+$\,.
\begin{lemma}\label{*AAnalyticity}
(i) There exist absolute constants $N_0\,,r_0\!>\!0$ such that for any $q\!\in\!\bH_{+\C}$ and
$n\!>\!N_0\|q\|_{\bH_{+\C}}$ the function $\psi_+(0,\cdot,q)$ has exactly $n$ roots, counted
with multiplicities, in the disc $\{\lambda:|\lambda|\!<\!4n\}$ and exactly one simple root in
the disc $\{\lambda:|\lambda\!-\!\sigma_n^0| \!<\! r_0n^{-1/2}\}$\,.

\noindent (ii) For each real potential $q\!\in\!\bH_+$ all eigenvalues $\sigma_n(q)$ extend
analytically to some complex ball $\{p\!\in\!\bH_{+\C}:\|p-q\|_{\bH_{+\C}}\!<\!R(q)\}$. Its
gradients\begin{footnote}{Recall that $\partial\xi(q)\big/\partial q\!=\!\zeta(q)$ means that
for any $v\!\in\!L^2$ the equation $(d_q\xi)(v)\!=\!(v,\overline{\zeta})_{L^2}$ holds true.
}\end{footnote} are given by
$$
\frac{\partial\sigma_n(q)}{\partial q(t)}= \psi_{n,D}^2(t,q)\,,
$$
where $\psi_{n,D}$ is the $n$-th normalized eigenfunction of $T_D$\,. Moreover,
$$
\frac{\partial \log[(-1)^n\psi'_+(0,\sigma_n(q),q)]}{\partial q(t)}=
-(\psi_{n,D}\chi_{n,D})(t,q)\,,
$$
where
$$
\chi_{n,D}(t,q)= \frac{\vartheta(t,\sigma_n(q),q)}{\psi'_{n,D}(0,q)}-
\frac{\dot{\psi}'_+}{\psi'_+}\,(0,\sigma_n(q),q)\cdot\psi_{n,D}(t,q)\,.
$$
\end{lemma}
\begin{proof} (i) The proof repeats the proof of \cite{CKK1} Lemma 4.1
(see also \cite{CK} Lemma 3.1).

\noindent (ii) The proof of the analyticity repeats the proof of \cite{CKK} Lemma 2.3 (p.
172). In order to calculate gradients note that the standard arguments (see \cite{CKK} Lemma
5.6) give
$$
\frac{\partial \psi_+(0,\lambda,q)}{\partial q(t)}=(\vp\p_+)(t,\lambda,q),\qquad
\frac{\partial \psi'_+(0,\lambda,q)}{\partial q(t)}=-(\vt\p_+)(t,\lambda,q),\quad t\!\ge\! 0.
$$
Applying the Implicit Function Theorem to the equation $\psi_+(0,\sigma_n(q),q)\!=\!0$ and
using the identity
$\int_{\R_+}{\psi_+^2}(t,\lambda,q)dt\!=\!\{\psi_+,\dot{\psi}_+\}(0,\lambda,q)$, we obtain
$$
\frac{\partial \sigma_n(q)}{\partial q(t)}=-\frac{\partial \psi_+(0)/\partial q(t)}{\partial
\psi_+(0)/\partial \lambda}= -\frac{(\varphi\psi_+)(t)}{\dot{\psi}_+(0)}=
-\frac{\psi_+^2(t)}{\psi_+'(0)\dot{\psi}_+(0)}=
\frac{\psi_+^2(t)}{\{\psi_+,\dot{\psi}_+\}(0)}=\psi_{2n+1}^2(t)
$$
and
$$
\frac{\partial\log[(-1)^n\psi'_+(0,\sigma_n(q),q)]}{\partial q(t)}=
\frac{-(\vartheta\psi_+)(t)+\dot{\psi}'_+(0)\cdot\psi_{2n+1}^2(t)}{\psi'_+(0)}=
-(\psi_{n,D}\chi_{n,D})(t)\,,
$$
where we omit $\sigma_n(q)$ and $q$ for short.
\end{proof}

\begin{lemma}\label{*AGradProd} (i) For each $q\in\bH_+$ and $n,m\!\ge\!0$ the following
identities are fulfilled:
$$
\begin{array}{ll}\displaystyle
\left((\psi_{n,D}^2)',\psi_{m,D}^2\right)_{+}=0\,,& \displaystyle
\left((\psi_{n,D}\chi_{n,D})',\psi_{m,D}^2\right)_{+}= {\textstyle -
\frac{1}{2}}\,\delta_{mn}\,, \vphantom{\big|_{\big|}}\cr \displaystyle
\left((\psi_{n,D}^2)',\psi_{m,D}\chi_{2m+1}\right)_{+}=
{\textstyle\frac{1}{2}}\,\delta_{mn}\,,& \displaystyle
\left((\psi_{n,D}\chi_{n,D})',\psi_{m,D}\chi_{m,D} \vphantom{\psi_{m,D}^2}\right)_{+}=0\,.
\vphantom{\big|^{\big|}}\end{array}
$$
\remark {\rm This Lemma is similar to \cite{CKK} Lemma 2.6 (see also \cite{PT} p. 44-45).}
\end{lemma}
\begin{proof}
For instance, we prove the third identity. Integration by parts gives
$$
I_{nm}=\int_{\R_+}(\psi_{n,D}^2)'(t,q)(\psi_{m,D}\chi_{m,D})(t,q)dt=
\frac{1}{2}\int_{\R_+}\{\psi_{m,D}\chi_{m,D}\,,(\psi_{n,D}^2)\}(t,q)dt
$$
$$
=\frac{1}{2}\int_{\R_+}(\chi_{m,D}\psi_{n,D}\{\psi_{m,D}\,,\psi_{n,D}\}+
\psi_{m,D}\psi_{n,D}\{\chi_{m,D}\,,\psi_{n,D}\})(t,q)dt\,.
$$
For $n\!\ne\!m$\,, this implies
$$
I_{nm}=\frac{1}{2(\sigma_m(q)\!-\!\sigma_n(q))}
(\{\psi_{m,D}\,,\psi_{n,D}\}\{\chi_{m,D}\,,\psi_{n,D}\})'(x,q,t)\Big|_{x=0}^{+\infty}=0\,.
$$
If $n\!=\!m$\,, then $\{\psi_{m,D}\,,\psi_{n,D}\}\!=\!0$\,,
$\{\chi_{m,D}\,,\psi_{n,D}\}\!=\!1$\,. Hence,
$I_{nn}\!=\!\frac{1}{2}\int_{\R_+}\psi_{n,D}^2(t,q)dt\!=\!\frac{1}{2}$\,.
\end{proof}

\noindent{\bf A.4 The leading terms of asymptotics of $\psi_+(0,\lambda)$ and
$\psi'_+(0,\lambda)$.}

\vskip 6pt

\noindent Recall that $\kappa_n\!=\!\p_{+}^0(0,\l_n^0)$\,,
$\kappa'_n\!=\!(\p_{+}^0)'(0,\l_n^0)$\,, $\dot{\kappa}_n\!=\!\dot{\p}_{+}^0(0,\l_n^0)$ and so
on. (\ref{Psi+00lambda}) yields
\begin{equation}
\label{KappaValues} \kappa_{2n+1}=0\,, \qquad
\kappa'_{2n+1}\asymp|\lambda_{2n+1}^0|^{1/4}\cdot a(\lambda_{2n+1}^0)\,, \qquad
\dot{\kappa}_{2n+1}\asymp|\lambda_{2n+1}^0|^{-1/4}\cdot a(\lambda_{2n+1}^0)\,.
\end{equation}

\begin{lemma}\label{*LKappaSpecAsymp}
The following asymptotics are fulfilled:
$$
\begin{array}{ll} \displaystyle
\frac{\dot{\kappa}'_{2n+1}}{\kappa'_{2n+1}}-
\frac{\ddot{\kappa}_{2n+1}}{2\dot{\kappa}_{2n+1}}=O(n^{-1})\,, & \displaystyle
\frac{\ddot{\kappa}'_{2n+1}}{\kappa'_{2n+1}}-
\frac{(\dot{\kappa}'_{2n+1})^2}{(\kappa'_{2n+1})^2} + \frac{\pi^2}{16} =O(n^{-1})\,,
\vphantom{\Big|_\Big|} \cr \vphantom{\Big|^\Big|} \displaystyle
\frac{\dddot\kappa_{2n+1}}{3\dot\kappa_{2n+1}}-
\frac{(\ddot\kappa_{2n+1})^2}{4(\dot\kappa_{2n+1})^2}+\frac{\pi^2}{48}=O(n^{-1}) &
\displaystyle \qquad as \quad n\!\to\!\infty\,.
\end{array}
$$
\end{lemma}
\begin{proof}
We prove the last asymptotics, the others are similar. Let
$f(\lambda)\!=\!\pi^{-\frac{1}{2}}2^{\frac{\lambda-1}{4}}\Gamma(\frac{\lambda+1}{4})$\,.
Identity (\ref{Psi+00lambda}) yields
$$
\dot{\kappa}_{2n+1}= \frac{(-1)^{n+1}\pi}{4}\,f(4n\!+\!3)\,,\qquad \ddot{\kappa}_{2n+1}=
2\cdot \frac{(-1)^{n+1}\pi}{4}\cdot \frac{df(\lambda)}{d\lambda}\Big|_{\lambda=4n+3}
$$
and
$$
\dddot{\kappa}_{2n+1}= \frac{(-1)^{n}\pi^3}{64}\,f(4n\!+\!3) + 3\cdot
\frac{(-1)^{n+1}\pi}{4}\cdot \frac{d^2f(\lambda)}{d\lambda^2}\Big|_{\lambda=4n+3}\,.
$$
Hence,
$$
\frac{\dddot\kappa_{2n+1}}{3\dot\kappa_{2n+1}}-
\frac{(\ddot\kappa_{2n+1})^2}{4(\dot\kappa_{2n+1})^2}+\frac{\pi^2}{48}=
\biggl[\frac{d^2f(\lambda)}{d\lambda^2}-\biggl(\frac{df(\lambda)}{d\lambda}\biggr)^2\biggr]
\Big|_{\lambda=4n+3}= \frac{d^2}{d\lambda^2}\,\log f(\lambda)\Big|_{\lambda=4n+3}=O(n^{-1})
$$
since $\frac{d^2}{dx^2}\log\Gamma(x)=O(x^{-1})$ as $x\!\to\!+\infty$\,.
\end{proof}

\begin{lemma} \label{*LA511}
Let $\psi_+^{(1)}\!=\!\psi_+^{(1)}(0,\sigma_n^0,q)$\,,
$\dot{\psi}_+^{(1)}\!=\!\dot{\psi}_+^{(1)}(0,\sigma_n^0,q)$ and so on. For some absolute
constant $\delta\!>\!0$ and all $q\!\in\!\bH_+$ the  following identities and asymptotics are
fulfilled:
$$
\begin{array}{ll}
\p_{+}^{(1)}=-2\dot{\kappa}_{2n+1}\nh q_{2n+1}\,, &
\dot{\p}_{+}^{(1)}\!=\!-2\dot{\kappa}_{2n+1}(\frac{\dot{\kappa}'_{2n+1}}{\kappa'_{2n+1}}\,\nh
q_{2n+1} \!+\! \frac{1}{2}\,\nc q_{2n+1} \!+ \el2_{\frac{1}{4}+\delta}(n)),
\vphantom{\big|_\big|}\cr (\p_{+}^{(1)})'=-\kappa'_{2n+1}\cdot \nc q_{2n+1}\,,
\vphantom{\Big|_\big|^\big|}&
(\dot{\p}_{+}^{(1)})'\!=\!\kappa'_{2n+1}(-\frac{\dot{\kappa}'_{2n+1}}{\kappa'_{2n+1}}\,\nc
q_{2n+1} \!+\! \frac{\pi^2}{8}\nh q_{2n+1} \!+ \el2_{\frac{1}{4}+\delta}(n)), \cr
\p_{+}^{(2)}\!=\!-2\dot{\kappa}_{2n+1}(-\nc q_{2n+1}\nh q_{2n+1} \!+
\el2_{\frac{3}{4}+\delta}(n)), & (\p_{+}^{(2)})' \!=\! \kappa'_{2n+1}(\frac{1}{2}(\nc
q_{2n+1})^2 \!-\! \frac{\pi^2}{8}(\nh q_{2n+1})^2 \!+ \el2_{\frac{3}{4}+\delta}(n)),
\vphantom{\big|^\big|}
\end{array}
$$
uniformly on bounded subsets of $\bH_+$\,.
\end{lemma}
\begin{proof} See \cite{CKK} Lemmas 5.11, 5.12 and \cite{CKK} Theorem 6.4.
\end{proof}

\noindent{\bf A.5 Three technical lemmas.}

\begin{lemma} \label{*xALemma}
Let $q,p\!\in\!\bH_+$\,. Then for all $\ve\!\in\!(0,\frac{1}{8})$ the following asymptotics is
fulfilled:
\begin{equation}
\label{xLA} a_{nm}=(q,(\psi_n\psi_m)(p))_+^2=\cases{O(n^{-\frac{1}{2}}\,m^{-\frac{1}{2}})& for
all $n,m\!\ge\!0$\,,\cr O(n^{-\frac{1}{2}-\frac{\varepsilon}{2}}m^{-\frac{1}{2}}), & if
$m\!\ge\!n\!+\!n^{\frac{1}{2}+\ve}$.}
\end{equation}
\end{lemma}

\begin{proof}
Using Lemma \ref{*AEstim}, Corollary \ref{*ADotEstim} and asymptotics (\ref{EstimExact}), we
obtain
$$
\psi_n(x,p)=\sqrt{2}\psi_n^0(x)+O(n^{-\frac{1}{2}}\log n\cdot\rho^{-1}(x,\lambda_n^0))\,,
\qquad \psi_n^0(x)=O(\rho^{-1}(x,\lambda_n^0))\,,
$$
where $\psi_n^0$ is the $n$-th unperturbed eigenfunction of the harmonic oscillator on $\R$.
Note that
$$
\int_0^{+\infty}\!\!\frac{dt}{(t^{1+\ve}\!+\!1)\cdot\rho^4(t,\lambda_n^0)}\,\le
\int_0^{+\infty}\!\!\frac{dt}{(t^{1+\ve}\!+\!1)(1\!+\!|\lambda_n^0\!-\!t^2|)}=O(n^{-1})\,.
$$
Therefore, for each $R\!\ge\!0$ we have
$$
\biggl[ \int_{R}^{+\infty}\!\!\!\!\frac{|q(t)|dt}{\rho(t,\lambda_n^0)\rho(t,\lambda_m^0)}
\biggr]^2 \le O(n^{-\frac{1}{2}}m^{-\frac{1}{2}})\cdot\int_R^{+\infty}\!\!
(t^{1+\ve}\!+\!1)\,|q(t)|^2dt =
O\left(n^{-\frac{1}{2}}m^{-\frac{1}{2}}(R\!+\!1)^{-(1-\ve)}\right)\,.
$$
If $R\!=\!0$\,, this gives $a_{nm}\!=\!O(n^{-\frac{1}{2}}m^{-\frac{1}{2}})$\,. Let
$m\!-\!n\!\ge\!n^{\frac{1}{2}+\ve}$. In this case we put $R\!=\!n^\ve$ and obtain
$$
a_{nm} = \left(q,(\psi_{n}^0\psi_{m}^0)\right)_{+}^2 + O(\log n\cdot n^{-1}m^{-\frac{1}{2}})
=\biggl[{\int_{0}^{n^\ve}}\!q(t)\psi_n^0(t)\psi_m^0(t)dt\biggr]^2+
O(n^{-\frac{1}{2}-(1-\ve)\ve}m^{-\frac{1}{2}})\,.
$$
Using WKB-bounds, it is easy to see (e.g., see \cite{CKK} Lemma 6.7) that
$$
\psi_n^0(t)= \sqrt{2\big/\pi}\cdot(\lambda_n^0)^{-\frac{1}{4}}\cos(\sqrt{\lambda_n^0}\cdot
t-{\textstyle\frac{\pi n}{2}}) + O(n^{-\frac{3}{4}+3\ve})\,,\qquad
|t|\!\le\!n^{\varepsilon}\,.
$$
Since $m\!\ge\!n\!+\!n^{\frac{1}{2}+\ve}$\,, we have
$\sqrt{\lambda_m^0}-\!\sqrt{\lambda_n^0}\ge n^{\ve}$\,. Integration by parts and
$q'\!\in\!L^2(\R)$ imply
$$
\int_{0}^{n^\ve}\!\!q(t)\psi_n^0(t)\psi_m^0(t)dt= O(n^{-\frac{1}{4}-\ve}m^{-\frac{1}{4}})\cdot
\int_{0}^{n^\ve}\!\!|q'(t)|dt + O(n^{-\frac{3}{4}+3\ve}m^{-\frac{1}{4}})=
O(n^{-\frac{1}{4}-\frac{\ve}{2}}m^{-\frac{1}{4}})\,.
$$
Therefore, $a_{nm}\!=\! O(n^{-\frac{1}{2}-\ve}m^{-\frac{1}{2}})\!+\!
O(n^{-\frac{1}{2}-(1-\ve)\ve}m^{-\frac{1}{2}})\!=\!
O(n^{-\frac{1}{2}-\frac{\ve}{2}}m^{-\frac{1}{2}})$\,.
\end{proof}

\begin{lemma}\label{*ALemma}
Let $\{a_{nm}\}_{n,m\ge 0}$ satisfy asymptotics (\ref{xLA}) for some
$\ve\!\in\!(0,\frac{1}{2})$\,. Then
$$
S_k= \sum_{n=0}^k\sum_{m=k+1}^{+\infty} \frac{a_{nm}}{m\!-\!n}\to 0\quad as \quad
k\!\to\!\infty\,.
$$
\end{lemma}
\begin{proof}
Let
$$
S_k^{(1)}\!=\!\sum_{n=0}^k\sum_{m=k+1}^{k+k^{\frac{1}{2}+\ve}}\frac{a_{nm}}{m-n}\qquad {\rm
and}\qquad
S_k^{(2)}\!=\!\sum_{n=0}^k\sum_{m=k+k^{\frac{1}{2}+\ve}}^{+\infty}\frac{a_{nm}}{m-n}\,.
$$
Using the simple estimate
$$
\sum_{m=k+1}^{k+k^{\frac{1}{2}+\ve}}\frac{1}{m\!-\!n}=
O\biggl(\log\frac{k\!+\!k^{\frac{1}{2}+\ve}\!-\!n}{k\!+\!1\!-\!n}\biggr)=
O\biggl(\frac{k^{\frac{1}{2}+\ve}}{k\!+\!1\!-\!n}\biggr)\,,
$$
we get
$$
S_k^{(1)}=
O(k^{\ve})\sum_{n=0}^k O(n^{-\frac{1}{2}}(k\!+\!1\!-\!n)^{-1})= O(k^{-\frac{1}{2}+\ve}\log
k)\,.
$$

Also, we have
$$
S_k^{(2)}=\sum_{n=0}^k\sum_{m=k+k^{\frac{1}{2}+\ve}}^{+\infty}
\frac{O(n^{-\frac{1}{2}-\frac{\ve}{2}}m^{-\frac{1}{2}})}{m\!-\!n}\le \sum_{n=0}^k
O(n^{-\frac{1}{2}-\frac{\ve}{2}})\sum_{m=k+1}^{+\infty}\frac{O(m^{-\frac{1}{2}})}{m\!-\!n}\,.
$$
Note that $\sum_{n=0}^k
O(n^{-\frac{1}{2}-\frac{\ve}{2}})\!=\!O(k^{\frac{1}{2}-\frac{\ve}{2}})$ and
$\sum_{m=k+1}^{+\infty}\frac{O(m^{-\frac{1}{2}})}{m-n}\!=\!O(\log k\cdot k^{-\frac{1}{2}})$\,.
Hence, \linebreak $S_k^{(2)}\!=\!O(k^{-\frac{\varepsilon}{2}}\log k)$ as $k\!\to\!\infty$\,.
Summarizing, we obtain $S_k\!\to\!0$ as $k\!\to\!\infty$\,.
\end{proof}

\begin{lemma}
\label{*ALaboutL2} Let $h_{n}\!=\!O(n^{-\beta})$ and $h_{n}\!=\!v\cdot (n\!+\!1)^{-\beta}+
\el2_{\beta-\frac{1}{4}}(n)$\,, where $v\!\in\!\R$ and $\beta\!\in\![\frac{1}{2}\,,1]$\,.
Then,
$$
\sum_{m:m\ne n} \frac{h_{m}}{(n\!-\!m)^2}=\frac{\pi^2v}{3(n\!+\!1)^\beta}+
\el2_{\beta-\frac{1}{4}}(n)\,.
$$
\end{lemma}
\begin{proof}
Let
$$ \sum_{m:m\ne n} \frac{1}{(m\!+\!1)^{\beta}(n\!-\!m)^{2}}=\sum_{m:|m-n|\le\sqrt{n}}\, +
\sum_{m:|m-n|>\sqrt{n}}= S_1+S_2\,.
$$
Since
$\sum_{m:|m-n|\le\sqrt{n}}\,\,(n\!-\!m)^{-2}=\frac{1}{3}\,\pi^2\!+\!O(n^{-\frac{1}{2}})$\,, we
have
$$
S_1=\biggl(\frac{1}{(n\!+\!1)^\beta}+O\left(n^{-\beta-\frac{1}{2}}\right)\biggr)\cdot
\biggl(\frac{\pi^2}{3}+O(n^{-\frac{1}{2}})\biggr)=
\frac{\pi^2}{3(n\!+\!1)^\beta}+O\left(n^{-\beta-\frac{1}{2}}\right)\,.
$$
Note that $S_2\!=\!O(n^{-1-\beta})$, if $\beta\!<\!1$, and $S_2\!=\!O(n^{-2}\log n)$, if
$\beta\!=\!1$\,. In any case, we obtain
$$
\sum_{m:m\ne n} \frac{v\cdot(m\!+\!1)^{-\beta}}{(n\!-\!m)^2}=
\frac{\pi^2v}{3(n\!+\!1)^\beta}+\el2_{\beta-\frac{1}{4}}(n)\,.
$$

Let $\widetilde h_m=h_m\!-\!v\cdot (m\!+\!1)^{-\beta}$.
We have
$$
\biggl(\sum_{m:m\ne n} \frac{\widetilde h_m} {(n\!-\!m)^2} \biggr)^{\!2} \le \sum_{m:m\ne n}
\frac{1}{(m\!+\!1)^{2\beta-\frac{1}{2}}(n\!-\!m)^2}\ \cdot\sum_{m:m\ne n}
\frac{(m\!+\!1)^{2\beta-\frac{1}{2}}\,\widetilde h_m^2}{(n\!-\!m)^2}\,.
$$
Using the simple estimate $\sum_{m:m\ne
n}{(m\!+\!1)^{-2\beta+\frac{1}{2}}(n\!-\!m)^{-2}}\!=\!O((n\!+\!1)^{-2\beta+\frac{1}{2}})$\,,
we deduce that
$$
\sum_{n\ge 0}\biggl((n\!+\!1)^{\beta-\frac{1}{4}}\!\!\!\sum_{m:m\ne n} \frac{\widetilde h_m}
{(n\!-\!m)^2} \biggr)^{\!2} \le O(1)\cdot \sum_{m\ge
0}(m\!+\!1)^{2\beta-\frac{1}{2}}\,\widetilde h_m^2\sum_{n:n\ne m}\frac{1}{(n\!-\!m)^2}= O(1)
$$
since $\widetilde h_m\!=\!\el2_{\beta-\frac{1}{4}}(m)$\,. Therefore, $\sum_{m:m\ne n}
{\widetilde h_m}\cdot {(n\!-\!m)^{-2}}=\el2_{\beta-\frac{1}{4}}(n)$\,.
\end{proof}

\noindent {\bf Acknowledgments.}\ Evgeny Korotyaev was partly supported by DFG project
BR691/23-1. Some part of this paper was written at the Mittag-Leffler Institute, Stockholm.
The authors are grateful to the Institute for the hospitality.

\end{document}